\documentclass[11pt]{article} 
\usepackage{theorem}
\usepackage{latexsym}
\usepackage{amssymb}
\usepackage{amsmath}
\usepackage{amsfonts}
\usepackage{amscd}
\usepackage{mathrsfs}
\usepackage{bezier}
\usepackage{color}
\usepackage{eufrak}
\usepackage{epic}
\makeatletter

\newbox\authrun
\newtoks\authorrunning
\newtoks\tocauthor
\newbox\titrun
\newtoks\titlerunning
\newtoks\toctitle

\def\institute#1{\gdef\@institute{#1}}
\def\author#1{\gdef\@author{#1}}
\def\maketitle{
\begin{center}
{\LARGE\bf \@title \par}
\vskip 2em
\end{center}
\let\maketitle=\relax}

\def\section{
\@startsection{section}
{1}{\z@}%{-7.5mm}
{3.5ex plus 1ex minus .2ex}
{2.3ex plus .2ex}
{\bf}}
\def\ps@run{%
 \def\@oddfoot{}\def\@evenfoot{}%       No feet.
 \def\@oddhead{\hfil \the\titlerunning/ \the\authorrunning\hfil \thepage }%
 \def\@evenhead{\hfil \hfil \thepage }}%

\def\footnoterule{\kern-3\p@
  \hrule width .4\columnwidth
  \kern 2.6\p@}
\long\def\@makefntext#1{\parindent 1em\noindent
            \hbox {\hss$\m@th^{\@thefnmark}$}#1}
            
\makeatother
\newcommand{\inv}{{\mathrm{inv}}}
\newcommand{\con}[2]{{\mathrm{con}_{#1}^{#2}}}
\newcommand{\res}[2]{{\mathrm{res}_{#1}^{#2}}}
\newcommand{\ind}[2]{{\mathrm{ind}_{#1}^{#2}}}

\newcommand{\set}{\mbox{-\boldmath $\mathrm{set}$}}%
\newcommand{\Map}{\mathrm{Map}}

\newcommand{\Mon}{\mbox{\boldmath$\mathrm{Mon}$}}
\newcommand{\El}{\mbox{\boldmath$\mathrm{El}$}}

\newcommand{\ACon}{\mbox{\boldmath $\mathrm{Con}$}_{\mathrm{alg}}}

\newcommand{\ARes}{\mbox{\boldmath $\mathrm{Res}$}_{\mathrm{alg}}}

\newcommand{\Green}{\mbox{\boldmath $\mathrm{Green}$}}
\newcommand{\obs}{{\mathrm{Obs}}}
\newcommand{\bz}{{\mathbb{Z}}}

\newcommand{\bq}{{\mathbb{Q}}}
\newcommand{\bc}{{\mathbb{C}}}

\newcommand{\id}{\mathrm{id}}
\renewcommand{\r}[1]{\,{}^{#1}}
\newcommand{\ov}[1]{\overline{#1}}
\newcommand{\st}{\mbox{\ }|\mbox{\ }}
\newcommand{\lr}[1]{\langle #1\rangle}
\renewcommand{\hom}{{\mathrm{Hom}}}
\renewcommand{\SS}{{\mathcal S}}
\newcommand{\SSS}{{\mathscr S}}
\newcommand{\C}{\mathrm{C}}

\renewcommand{\L}{{\mathscr L}}
\newcommand{\U}{{\mathcal U}}
\newcommand{\R}{{\mathcal R}}
\newcommand{\im}{{\mathrm{Im\hspace{0.03cm}}}}
\renewcommand{\ker}{{\mathrm{Ker}\hspace{0.05cm}}}
\newcommand{\XX}{{\mathfrak X}}

\newcommand{\coc}{Z^1}
\newcommand{\coh}{H^1}
\newcommand{\boldGamma}{\text{$\boldsymbol{\varGamma}$}}
\newcommand{\aut}{{\mathrm{Aut}}}
\newcommand{\chara}{{\mathrm{char}}}
\newcommand{\trans}{\mathbf{trans}}

\newcommand{\Hom}{{\mathrm{Hom}}}
\newtheorem{theorem}{Theorem}[section]
\newtheorem{proposition}[theorem]{Proposition}
\newtheorem{lemma}[theorem]{Lemma}
\newtheorem{corollary}[theorem]{Corollary}
\theorembodyfont{\rmfamily}
\newtheorem{definition}[theorem]{Definition}
\newtheorem{example}[theorem]{Example}

\newtheorem{condition}{(A)}

\newtheorem{remark}[theorem]{\mdseries{\itshape{Remark}}}
\newtheorem{notation}{\mdseries{\itshape{Notation}}}% and terminology
\def\theenumi{\alph{enumi}}
\numberwithin{equation}{section}
\begin{document}
\pagestyle{run}
\thispagestyle{plain}
%-------------------------------------------------------
\title{The lattice Burnside rings} 
%\author{Yugen Takegahara}
\institute{Muroran Institute of Technology, 27-1 Mizumoto, Muroran 050-8585, 
Japan \\ \rm E-mail: yugen@mmm.muroran-it.ac.jp}
\titlerunning{The lattice Burnside rings}
\authorrunning{F. Oda, Y. Takegahara, and T. Yoshida}
\maketitle
\vskip 1em
\centerline{\Large\bf Fumihito Oda}
\vskip 1em
\centerline{\small\itshape\lineskip .75em 
Department of Mathematics, Kindai University,
Higashi-Osaka, 577-8502, Japan}
\centerline{\small\lineskip .75em 
{\itshape E-mail:} odaf@math.kindai.ac.jp}
\vskip 2em
\centerline{\Large\bf Yugen Takegahara\!\,\footnotemark}
\vskip 1em
\centerline{\small\itshape\lineskip .75em 
Muroran Institute of Technology, 27-1 Mizumoto, Muroran 050-8585, Japan} 
\centerline{\small\lineskip .75em 
{\itshape E-mail:} yugen@mmm.muroran-it.ac.jp}
\vskip 2em
\centerline{\Large\bf Tomoyuki Yoshida}
\vskip 1em
\centerline{\small\lineskip .75em 
{\itshape E-mail:} ytomoyuki@mub.biglobe.ne.jp}
\vskip 2em
\begin{abstract}
We introduce the concept of lattice Burnside ring for a finite group $G$ 
associated to a family of nonempty sublattices $\L_H$ of a finite 
$G$-lattice $\L$ for $H\leq G$. The slice Burnside ring introduced by 
S. Bouc \cite{Bou12} is isomorphic to a lattice Burnside ring. Any lattice 
Burnside ring is isomorphic to an abstract Burnside ring. The ring structure 
of a lattice Burnside ring is explored on the basis of the fundamental theorem 
of abstract Burnside rings. We study the units and the primitive idempotents 
of a lattice Burnside ring. There are certain rings called partial 
lattice Burnside rings. Any partial lattice Burnside ring, which is isomorphic 
to an abstract Burnside ring, consists of elements of a lattice 
Burnside ring; it is not necessarily a subring. The section Burnside ring 
introduced by S. Bouc \cite{Bou12} is isomorphic to a partial lattice Burnside 
ring. 
\end{abstract}
%
%\tableofcontents
%
%footnote----------------------------------------------------------------
%
\footnotetext{
{\rm This work was supported by JSPS KAKENHI Grant Number JP16K05052.}
}
%----------------------------------------------------------------
%
%introduction
%
%----------------------------------------------------------------
\section{Introduction}
%----------------------------------------------------------------
Let $G$ be a finite group, and let $\L$ be a finite $G$-lattice, that is, 
$\L$ is a finite lattice on which $G$ acts and the binary relation in $\L$ 
is invariant under the action of $G$. The purpose of this paper is to 
introduce the concept of lattice Burnside ring associated to a family of 
nonempty sublattices $\L_H$ of $\L$ for $H\leq G$, together with its ring 
structure. The slice Burnside ring of $G$ introduced by 
S. Bouc \cite{Bou12}, which arises  from morphisms of finite $G$-sets, 
inspired us to study lattice Burnside rings. 
\par
In Section 2, we first introduce the concept of monoid functor 
$M=(M,\con{}{},\res{}{})$ assigning each $H\leq G$ a monoid $M(H)$. Let 
$M=(M,\con{}{},\res{}{})$ be a monoid functor, and let $H\leq G$. There is 
an additive contravariant functor $T^M_H$ from the category of finite left 
$H$-sets to the category of monoids such that for any finite left $H$-set $J$, 
$T^M_H(J)$ is the set of $H$-equivariant maps $\pi:J\to\dot\cup_{K\leq H}M(K)$ 
with $\pi(x)\in M(H_x)$, where $H_x$ is the stabilizer of $x$ in $H$. A pair 
$(J,\pi)$ of a left $H$-set $J$ and $\pi\in T^M_H(J)$ is called an element of 
$T^M_H$. The $M$-Burnside ring $\mathrm{\Omega}(H,M)$ is defined to be 
the Grothendieck ring of the category of elements of $T^M_H$ (cf. 
\cite{Jacob86,Na08,OY01}). 
\par
Any finite lattice is considered as a monoid with the binary operation 
given by `meet', so that $\L$ is viewed as a finite monoid. We introduce 
the concept of a monoid functor $M_\L=(M_\L,\con{}{},\res{}{})$ assigning each 
$H\leq G$ a nonempty sublattice $\L_H$ of $\L$ (see 
Proposition \ref{pro:lattice}), and call $\mathrm{\Omega}(G,M_\L)$ the lattice 
Burnside ring associated to the family of nonempty sublattices $\L_H$ of $\L$ 
for $H\leq G$. If $\sup\L$ is $G$-invariant, then there is a monoid functor 
$C_\L=(C_\L,\con{}{},\res{}{})$ assigning each $H\leq G$ the set of 
$H$-invariants in $\L$, and the lattice Burnside ring 
$\mathrm{\Omega}(G,C_\L)$ is isomorphic to the crossed Burnside ring 
associated to $\L$. We consider the set $\SSS:=\SSS(G)$ of subgroups of $G$ to 
be a finite $G$-lattice with the binary relation given by inclusion and 
the action of $G$ given by conjugation. There is a monoid functor 
$M_\SSS=(M_\SSS,\con{}{},\res{}{})$ assigning each $H\leq G$ the set of 
subgroups of $G$ containing $H$, and the lattice Burnside ring 
$\mathrm{\Omega}(G,M_\SSS)$ is isomorphic to the slice Burnside ring of $G$. 
We also obtain a monoid functor 
$M_\SSS^\circ=(M_\SSS^\circ,\con{}{},\res{}{})$ assigning each $H\leq G$ 
the set of normal subgroups of $H$, and direct our attention to the deference 
between $\mathrm{\Omega}(G,M_\SSS^\circ)$ and $\mathrm{\Omega}(G,M_\SSS)$. 
\par
The lattice Burnside rings are extensions of the ordinary Burnside ring of 
$G$. There are various attempts to generalize the theory of ordinary Burnside 
ring of $G$. Among others, the concept of abstract Burnside ring was 
introduced in \cite{YOT18}, and any lattice Burnside ring is isomorphic to 
an abstract Burnside ring. By using the fundamental theorem of abstract 
Burnside rings, we determine the primitive idempotents of the $\bq$-algebra 
$\bq\otimes_\bz\mathrm{\Omega}(G,M_\L)$ in Section 4, and establish 
a criterion for the units of $\mathrm{\Omega}(G,M_\L)$ which generalize 
\cite[Proposition 6.5]{Yo90a} in Section 5. 
\par
The theory of abstract Burnside rings provides a principle of determining 
the primitive idempotents of the lattice Burnside rings. In Section 6, we 
present a characterization of solvable groups for a certain class of lattice 
Burnside rings as a sequel to the study of primitive idempotents of 
$\mathrm{\Omega}(G,M_\L)$, which extend 
Dress' characterization of solvable groups for the ordinary Burnside ring of 
$G$ (cf. \cite{Dr69}). 
\par
In Section 7, we introduce the concept of partial lattice Burnside ring on 
the basis of results shown in \cite{Yo90b}, together with its ring structure. 
Any partial lattice Burnside ring, which is isomorphic to an abstract Burnside 
ring, consists of elements of a lattice Burnside ring; it is not necessarily 
a subring. The section Burnside ring of $G$ introduced by 
S. Bouc \cite{Bou12} is isomorphic to a partial lattice Burnside ring. 
\begin{notation}
Let $G$ be a finite group. The rational integers, the rational numbers, and 
the complex numbers are denoted by $\bz$, $\bq$, and $\bc$, respectively. We 
denote by $\epsilon$ the identity of $G$. The subgroup generated by an element 
$g$ of $G$ is denoted by $\lr{g}$. For subgroups $H$ and $K$, we write 
$K\leq H$ if $K$ is a subgroup of $H$. Let $H\leq G$. We denote by $H\set$ 
the category of finite left $H$-sets and $H$-equivariant maps. Given 
$J\in G\set$ and $x\in J$, $H_x$ denotes the stabilizer of $x$ in $H$. We set 
$\r{r}\!H=rHr^{-1}$ and $\r{r}\!g=rgr^{-1}$ for all $g,\,r\in G$. For each 
$K\leq H$, $H/K$ denotes the set of left cosets $hK$, $h\in H$, of $K$ in $H$. 
Given a pair $(U,\,K)$ of subgroups $U$ and $K$ of $H$, we denote by 
$U\backslash H/K$ the set of $(U,K)$-double cosets $UhK$, $h\in H$, in $H$. 
The identity map on a set $\Sigma$ is denoted by $\id_\Sigma$. For a finite 
set $\Sigma$, $|\Sigma|$ denotes the cardinality $\sharp\Sigma$ of $\Sigma$. 
For each positive integer $m$, we set $[m]=\{1,\,2,\dots,\,m\}$. 
\end{notation}
%----------------------------------------------------------------
%
%newsection
%
%----------------------------------------------------------------
\section{$M$-Burnside ring functors}
%----------------------------------------------------------------
We first review the Mackey functors and the Green functors (cf. 
\cite{Bol98,Gre71,Theve88,Yo80}). Let $G$ be a finite group, and let $k$ be 
a commutative unital ring. A Mackey functor for $G$ over $k$ is a quadruple 
$X=(X,\con{}{},\res{}{},\ind{}{})$ consisting of a family of $k$-modules 
$X(H)$ for $H\leq G$ and a family of $k$-module homomorphisms 
\[
\begin{array}{ll}
\vspace{0.1cm}
\con{H}{g}:X(H)\to X(\r{g}\!H), \quad&\mbox{the conjugation maps}, \\
\vspace{0.1cm}
\res{K}{H}:X(H)\to X(K), \quad&\mbox{the restriction maps}, \\
\ind{K}{H}:X(K)\to X(H), \quad&\mbox{the induction maps}, 
\end{array}
\]
for $g\in G$ and $K\leq H\leq G$, satisfying the axioms 
\[
\begin{array}{lll}
\vspace{0.1cm}
(\mathrm{G}.1)&
\con{\r{r}\!H}{g}\circ\con{H}{r}=\con{H}{gr}, \quad\con{H}{h}=\id_{X(H)},\\
\vspace{0.1cm}
(\mathrm{G}.2)
&\res{L}{K}\circ\res{K}{H}=\res{L}{H}, \quad\res{H}{H}=\id_{X(H)},\\
\vspace{0.1cm}
(\mathrm{G}.3)
&\con{K}{g}\circ\res{K}{H}=\res{\r{g}\!K}{\r{g}\!H}\circ\con{H}{g}, \\
\vspace{0.1cm}
(\mathrm{G}.4)
&\ind{K}{H}\circ\ind{L}{K}=\ind{L}{H}, \quad\ind{H}{H}=\id_{X(H)},\\
\vspace{0.1cm}
(\mathrm{G}.5)
&\con{H}{g}\circ\ind{K}{H}=\ind{\r{g}\!K}{\r{g}\!H}\circ\con{K}{g},\\
\vspace{0.1cm}
(\mathrm{G}.6)&\mbox{(Mackey axiom)}\\
\mbox{}&\displaystyle
\res{K}{H}\circ\ind{U}{H}=\sum_{KhU\in K\backslash H/U}
\ind{K\cap\r{h}\!U}{K}\circ\res{K\cap\r{h}\!U}{\r{h}\!U}\circ\con{U}{h}
\end{array}
\]
for all $U\leq H\leq G$, $L\leq K\leq H$, $g,\,r\in G$, and $h\in H$, which 
was introduced by Green \cite{Gre71}. A Green functor for $G$ over $k$ is 
a Mackey functor $X=(X,\con{}{},\res{}{},\ind{}{})$ for $G$ over $k$ such 
that the $k$-modules $X(H)$ for $H\leq G$ are $k$-algebras, the conjugation 
maps and the restriction maps are $k$-algebra homomorphisms, and the axioms 
\[
\begin{array}{ll}
\vspace{0.1cm}
(\mathrm{G}.7)&\mbox{(Frobenius axioms)}\\
\mbox{}&
x\cdot\ind{K}{H}(y)=\ind{K}{H}(\res{K}{H}(x)\cdot y), \quad
\ind{K}{H}(y)\cdot x=\ind{K}{H}(y\cdot\res{K}{H}(x)) 
\end{array}
\]
are satisfied for all $K\leq H$, $x\in X(H)$, and $y\in X(K)$. 
\par
A semigroup with identity is called a monoid. We denote by $\epsilon_A$ 
the identity of a monoid $A$. A map $f:A\to B$ between monoids is called 
a (monoid) homomorphism if $f(\epsilon_A)=\epsilon_B$ and 
$f(x\cdot y)=f(x)\cdot f(y)$ for all $x,\,y\in A$. 
\begin{definition}
\label{def:monoid}
A monoid functor for $G$ is a triple $M=(M,\con{}{},\res{}{})$ consisting 
of a family of monoids $M(H)$ for $H\leq G$ and a family of monoid 
homomorphisms 
\[
\begin{array}{ll}
\vspace{0.1cm}
\con{H}{g}:M(H)\to M(\r{g}\!H), \quad&\mbox{the conjugation maps}, \\
\res{K}{H}:M(H)\to M(K), \quad&\mbox{the restriction maps}, 
\end{array}
\]
for $g\in G$ and $K\leq H\leq G$, satisfying the axioms 
\[
\begin{array}{lll}
\vspace{0.1cm}
(\mathrm{M}.1)&
\con{\r{r}\!H}{g}\circ\con{H}{r}=\con{H}{gr}, \quad\con{H}{h}=\id_{M(H)}, \\
\vspace{0.1cm}
(\mathrm{M}.2)
&\res{L}{K}\circ\res{K}{H}=\res{L}{H}, \quad\res{H}{H}=\id_{M(H)},\\
\vspace{0.1cm}
(\mathrm{M}.3)
&\con{K}{g}\circ\res{K}{H}=\res{\r{g}\!K}{\r{g}\!H}\circ\con{H}{g} 
\end{array}
\]
for all $L\leq K\leq H\leq G$, $g,\,r\in G$, and $h\in H$. 
\end{definition}
\par
Henceforth, we suppose that $M$ is a monoid functor for $G$. Let $H\leq G$, 
and set 
\[
\widetilde{M}(H)=\dot{\displaystyle\bigcup_{K\leq H}}M(K). 
\]
We consider $\widetilde{M}(H)$ to be a left $H$-set with the action given 
\[
\r{h}\!s=\con{K}{h}(s)
\]
for all $h\in H$ and $s\in M(K)$ with $K\leq H$. Given $J,\,J_0\in H\set$, let 
$\Map_H(J_0,J)$ be the set of $H$-equivariant maps from $J_0$ to $J$. There 
exists an additive contravariant functor $T^M_H:H\set\to\Mon$, where $\Mon$ is 
the category of monoids, such that for each $J\in H\set$, $T^M_H(J)$ is 
defined to be the monoid 
\[
\{\pi\in\Map_H(J,\widetilde{M}(H))\st\pi(x)\in M(H_x)\quad\mbox{for all}\quad 
x\in J\} 
\]
with pointwise multiplication, where $H_x$ is the stabilizer of $x$ in $H$, 
and the morphism $T^M_H(\lambda):T^M_H(J)\to T^M_H(J_0)$ with 
$J,\,J_0\in H\set$ and $\lambda\in\Map_H(J_0,J)$ is defined by 
\[
T^M_H(\lambda)(\pi):J_0\to \widetilde{M}(H),\quad x\mapsto
\res{H_x}{H_{\lambda(x)}}(\pi(\lambda(x)))
\]
for all $\pi\in T^M_H(J)$ (cf. \cite[\S2]{Jacob86}). We call a pair 
$(J,\pi)$ of $J\in H\set$ and $\pi\in T^M_H(J)$ an element of $T^M_H$. 
The morphisms $\lambda:(J_0,\pi_0)\to(J,\pi)$ between elements 
$(J_0,\pi_0)$ and $(J,\pi)$ of $T^M_H$ 
are defined to be the $H$-equivariant maps $\lambda:J_0\to J$ such that 
$T^M_H(\lambda)(\pi)=\pi_0$. Thus we obtain the category of elements of 
$T^M_H$ (cf. \cite[(2.10)]{OY01}). 
\begin{definition}
\label{def:element}
Let $H\leq G$. For each element $(J,\pi)$ of $T^M_H$, let $[J,\pi]$ denote 
the isomorphism class of elements of $T^M_H$ containing $(J,\pi)$. We define 
a ring $\mathrm{\Omega}(H,M)$ to be the ring consisting of all 
$\bz$-linear combinations of isomorphism classes of elements of $T^M_H$ with 
addition and multiplication given by 
\[
[J_1,\pi_1]+[J_2,\pi_2]=[J_1\dot\cup J_2,\pi_1\dot+\pi_2]\quad\mbox{and}\quad
[J_1,\pi_1]\cdot[J_2,\pi_2]=[J_1\times J_2,\pi_1\cdot\pi_2], 
\]
where
\[
\pi_1\dot+\pi_2:J_1\dot\cup J_2\to \widetilde{M}(H),\quad
x\mapsto\pi_1(x)\;\mbox{ if }x\in J_1, \quad
x\mapsto\pi_2(x)\;\mbox{ if }x\in J_2
\]
and
\[
\pi_1\cdot\pi_2:J_1\times J_2\to \widetilde{M}(H),\quad(x_1,x_2)\mapsto
\pi_1(x_1)\cdot\pi_2(x_2), 
\]
where 
\[
\pi_1(x_1)\cdot\pi_2(x_2)=\res{H_{x_1}\cap H_{x_2}}
{H_{x_1}}(\pi_1(x_1))\cdot\res{H_{x_1}\cap H_{x_2}}{H_{x_2}}(\pi_2(x_2)), 
\]
for all elements $(J_1,\pi_1)$ and $(J_2,\pi_2)$ of $T^M_H$, and call it 
the $M$-Burnside ring. This ring is the $F$-Burnside ring with $F=T^M_H$ 
introduced by Jacobson \cite{Jacob86} (cf. \cite{Na08}). 
\end{definition}
\par
For each $H\leq G$, let $\bz M(H)$ be the $\bz$-algebra consisting of all 
$\bz$-linear combinations of elements of $M(H)$ with multiplication given by 
the binary operation of $M(H)$, which is called a monoid ring. Then 
the family of monoid rings $\bz M(H)$ for $H\leq G$ defines an algebra 
restriction functor $\bz M=(\bz M,\con{}{},\res{}{})$ with conjugation maps 
and restriction maps given by those of the monoid functor 
$M=(M,\con{}{},\res{}{})$ (cf. \cite{Bol98,TY14}). The ring 
$\mathrm{\Omega}(H,T^{\bz M}_H)$ defined in \cite[\S3]{TY14} coincides 
with $\mathrm{\Omega}(H,M)$. 
\par
Let $H\leq G$. The Burnside ring $\mathrm{\Omega}(H)$ is the commutative ring 
consisting of all $\bz$-linear combinations of isomorphism classes $[J]$ for 
$J\in H\set$ with disjoint union for addition and cartesian product for 
multiplication (cf. \cite[\S80]{CR87}). When $\widetilde{M}(H)$ is 
the monoid $\{\epsilon\}$ consisting of only the identity $\epsilon$, we 
regard $\mathrm{\Omega}(H,M)$ with $\mathrm{\Omega}(H)$. 
\par
The Burnside ring functor 
$\mathrm{\Omega}=(\mathrm{\Omega},\con{}{},\res{}{},\ind{}{})$, which is 
a Green functor, is defined to be the family of $\bz$-algebras 
$\mathrm{\Omega}(H)$ for $H\leq G$, together with conjugation maps 
$\con{H}{g}:\mathrm{\Omega}(H)\to\mathrm{\Omega}(\r{g}\!H)$, where $g\in G$, 
restriction maps $\res{K}{H}:\mathrm{\Omega}(H)\to\mathrm{\Omega}(K)$, and 
induction maps $\ind{K}{H}:\mathrm{\Omega}(K)\to\mathrm{\Omega}(H)$, where 
$K\leq H$ in both cases, arising from the usual conjugation $\con{H}{g}(J)$, 
restriction $\res{K}{H}(J)$, and induction $\ind{K}{H}(V)$, where $J\in H\set$ 
and $V\in K\set$ (cf. \cite[Example 2.11]{Yo80}). 
\par
Following \cite[\S3]{TY14}, we define the $M$-Burnside ring functor. 
Let $g\in G$, and let $K\leq H\leq G$. Given an element $(J,\pi)$ of $T^M_H$ 
and an element $(V,\varpi)$ of $T^M_K$, we define 
$\r{g}\!\pi\in T^M_{\r{g}\!H}(\con{H}{g}(J))$, 
$\pi|_K\in T^M_K(\res{K}{H}(J))$, and $\varpi^H\in T^M_H(\ind{K}{H}(V))$ by 
\[
(\r{g}\!\pi)(g\otimes x)=\r{g}\!\pi(x),\quad
\pi|_K(x)=\res{K_x}{H_x}(\pi(x)),\quad\mbox{and}\quad
\varpi^H(h\otimes y)=\r{h}\!\varpi(y)
\]
for all $x\in J$, $y\in V$, and $h\in H$ (see \cite[p. 97]{TY14}). 
\begin{definition}
\label{def:Burn-element}
For each $H\leq G$, let $\mathrm{\Omega}^M(H)$ denote the $M$-Burnside ring 
$\mathrm{\Omega}(H,M)$. We define a Green functor 
$\mathrm{\Omega}^M=(\mathrm{\Omega}^M,\con{}{},\res{}{},\ind{}{})$ to be 
the family of $\bz$-algebras $\mathrm{\Omega}^M(H)$ for $H\leq G$, together 
with conjugation maps 
$\con{H}{g}:\mathrm{\Omega}^M(H)\to\mathrm{\Omega}^M(\r{g}\!H)$, where 
$g\in G$, restriction maps 
$\res{K}{H}:\mathrm{\Omega}^M(H)\to\mathrm{\Omega}^M(K)$, and induction maps 
$\ind{K}{H}:\mathrm{\Omega}^M(K)\to\mathrm{\Omega}^M(H)$, where $K\leq H$ in 
both cases, given by 
\[
\begin{array}{l}
\vspace{0.3cm}
\con{H}{g}([J,\pi])=[\con{H}{g}(J),\r{g}\!\pi],\\
\vspace{0.3cm}
\res{K}{H}([J,\pi])=[\res{K}{H}(J),\pi|_K],\\
\ind{K}{H}([V,\varpi])=[\ind{K}{H}(V),\varpi^H] 
\end{array}
\]
for all elements $(J,\pi)$ of $T^M_H$ and for all elements $(V,\varpi)$ of 
$T^M_K$. 
\end{definition}
\par
We call the Green functor 
$\mathrm{\Omega}^M=(\mathrm{\Omega}^M,\con{}{},\res{}{},\ind{}{})$ 
the $M$-Burnside ring functor, which is a $G$-functor version of 
the $F$-Burnside ring functor defined in \cite{Jacob86,Na08}. 
\par
Let $H\leq G$, and set $\SS(H,M)=\{(K,s)\st K\leq H,\,s\in M(K)\}$. For each 
$(K,s)\in\SS(H,M)$, we define an $H$-equivariant map 
$\pi_s:H/K\to\widetilde{M}(H)$ by 
\[
hK\mapsto\r{h}\!s\in M(\r{h}\!K) 
\]
for all $h\in H$, and set $[(H/K)_s]=[H/K,\pi_s]\in\mathrm{\Omega}(H,M)$. 
\begin{definition}
\label{def:pi}
Given $H\leq G$, we define an action of $H$ on $\SS(H,M)$ by 
\[
h.(K,s)=(\r{h}\!K,\r{h}\!s) 
\]
for all $h\in H$ and $(K,s)\in\SS(H,M)$, and denote by $\R(H,M)$ a complete 
set of representatives of $H$-orbits in $\SS(H,M)$ such that $K\in\C(H)$, 
where $\C(H)$ is a full set of nonconjugate subgroups of $H$, for each 
$(K,s)\in\R(H,M)$. 
\end{definition}
\par
Given $(K,s),\,(U,t)\in\SS(H,M)$, it is easily verified that 
$(H/K)_s\simeq(H/U)_t$ if and only if $(K,s)=h.(U,t)$ for some $h\in H$ 
(cf. Example \ref{epl:trans}). 
\begin{proposition}
\label{pro:monoid}
Let $H\leq G$. The elements $[(H/K)_s]$ for $(K,s)\in\R(H,M)$ form 
a $\bz$-basis of $\mathrm{\Omega}(H,M)$, and multiplication on 
$\mathrm{\Omega}(H,M)$ is given by 
\[
[(H/K)_s]\cdot[(H/U)_t]
=\sum_{KhU\in K\backslash H/U}\left[(H/(K\cap\r{h}U))_
{\res{K\cap\r{h}U}{K}(s)\cdot\res{K\cap\r{h}U}{\r{h}U}(\r{h}t)}\right] 
\]
for all $(K,s),\,(U,t)\in\R(H,M)$. Moreover, 
\[
\begin{array}{l}
\vspace{0.3cm}
\con{H}{g}([(H/U)_t])=[(\r{g}\!H/\r{g}U)_{\r{g}t}],\\
\vspace{0.3cm}
\displaystyle\res{K}{H}([(H/U)_t])
=\sum_{KhU\in K\backslash H/U}\left[(K/(K\cap\r{h}U))_
{\res{K\cap\r{h}U}{\r{h}U}(\r{h}t)}\right],\\
\ind{K}{H}([(K/L)_s])=[(H/L)_s] 
\end{array}
\]
for all $g\in G$, $(U,t)\in\R(H,M)$, $K\leq H$, and $(L,s)\in\R(K,M)$. 
\end{proposition}
{\itshape Proof.}
The proof is straightforward. (See also the proof of 
\cite[Proposition 3.1]{TY14}.) 
$\Box$
\begin{example}
\label{epl:monomial}
Let $A$ be a finite abelian $G$-group, that is, a finite abelian group on which $G$ acts as automorphisms of $A$ 
(cf. \cite[Chapter 1, Definition 8.1]{Su82}). Given $g\in G$ and $a\in A$, 
the effect of $g$ on $a$ is denoted by $\r{g}\!a$. Let $H\leq G$. By 
restriction of operators from $G$ to $H$, we view $A$ as an $H$-group. A map 
$\sigma:H\to A$ is called a $1$-cocycle or a crossed homomorphism if 
\[
\sigma(h_1h_2)=\sigma(h_1)\r{h_1}\!\sigma(h_2)
\]
for all $h_1,\,h_2\in H$ (cf. \cite[I, p. 243]{Su82}). The set $\coc(H,A)$ of 
$1$-cocycles from $H$ to $A$ is an abelian group with the product operation 
given by 
\[
\sigma\cdot\tau(h)=\sigma(h)\tau(h)
\]
for all $\sigma,\,\tau\in\coc(H,A)$ and $h\in H$. Let $\sigma\in\coc(H,A)$. 
For each $K\leq H$, let $\sigma|_K:K\to A$ denote the $1$-cocycle obtained by 
restriction of $\sigma:H\to A$ from $H$ to $K$. Given $g\in G$ and $a\in A$, 
we define $1$-cocycles $g\sigma:\r{g}\!H\to A$ and $\sigma^a:H\to A$ by 
\[
(g\sigma)(ghg^{-1})=\r{g}\!\sigma(h)\quad\mbox{and}\quad
\sigma^a(h)=a^{-1}\sigma(h)\r{h}\!a
\]
for all $h\in H$, respectively. Let $\ov{\sigma}$ denote the $A$-orbit 
$\{\sigma^a\st a\in A\}$ in $\coc(H,A)$ containing $\sigma$. We denote by 
$\coh(H,A)$ the set of $A$-orbits in $\coc(H,A)$, that is, 
\[
\coh(H,A)=\{\ov{\sigma}\st\sigma\in\coc(H,A)\}, 
\]
and make it into an abelian group by defining 
\[
\ov{\sigma}\cdot\ov{\tau}=\ov{\sigma\cdot\tau} 
\]
for all $\sigma,\,\tau\in\coc(H,A)$. Observe that 
$\ov{g(\sigma^a)}=\ov{(g\sigma)^{\r{g}\!a}}=\ov{g\sigma}$ for 
all $g\in G$ and $a\in A$. We define a monoid functor 
$H_A=(H_A,\con{}{},\res{}{})$ for $G$ by 
\[
\begin{array}{l}
\vspace{0.3cm}
H_A(H)=\coh(H,A), \\
\vspace{0.3cm}
\con{H}{g}:\coh(H,A)\to\coh(\r{g}\!H,A),\;\ov{\sigma}\mapsto\ov{g\sigma},\\
\res{K}{H}:\coh(H,A)\to\coh(K,A),\;\ov{\sigma}\mapsto\ov{\sigma|_K}
\end{array}
\]
for all $g\in G$ and $K\leq H\leq G$. The $H_A$-Burnside ring 
$\mathrm{\Omega}(G,H_A)$ is isomorphic to the ring of monomial representations 
of $G$ with coefficients in $A$ introduced by Dress \cite{Dr71c} which is also 
called the monomial Burnside ring for $G$ with fibre group $A$ in 
\cite{Bark04}. 
\end{example}
\begin{example}
\label{epl:crossed}
Let $S$ be a finite $G$-monoid, that is, a finite monoid on which $G$ acts as 
monoid homomorphisms (cf. \cite[(2.1)]{OY01}). Given $g\in G$ and $s\in S$, 
$\r{g}\!s$ denotes the effect of $g$ on $s$. We define a monoid functor 
$C_S=(C_S,\con{}{},\res{}{})$ for $G$ by 
\[
\begin{array}{l}
\vspace{0.3cm}
C_S(H)=\{s\in S\st\r{h}\!s=s\quad\mbox{for all}\quad h\in H\}, \\
\vspace{0.3cm}
\con{H}{g}:C_S(H)\to C_S(\r{g}\!H),\;s\mapsto\r{g}\!s,\\
\res{K}{H}:C_S(H)\to C_S(K),\;s\mapsto s 
\end{array}
\]
for all $g\in G$ and $K\leq H\leq G$. The crossed Burnside ring functor 
defined in \cite[\S4]{OY04} is isomorphic to the $C_S$-Burnside ring functor . 
The ring $\mathrm{\Omega}(G,C_S)$ is isomorphic to the crossed Burnside ring 
of $G$ associated to $S$ (cf. \cite{Bou03b, OY01}). 
\end{example}
%----------------------------------------------------------------
%
%newsection
%
%----------------------------------------------------------------
\section{A fundamental theorem of $M$-Burnside rings}
%----------------------------------------------------------------
Following \cite[\S4]{OY01}, we present a fundamental theorem for 
$M$-Burnside rings. The results in this section are special cases of those in 
\cite[\S9]{TY14}. 
\par
Let $H\leq G$. For any $g\in G$, there is a map 
$\con{H}{g}:\bz M(H)\to\bz M(\r{g}\!H)$ given by 
\[
\sum_j\ell_js_j\mapsto\sum_j\ell_j\r{g}\!s_j 
\]
for all $\sum_j\ell_js_j\in\bz M(H)$ with $\ell_j\in\bz$ and $s_j\in M(H)$. 
\par
We define a subring $\mho(H,M)$ of $\prod_{U\leq H}\bz M(U)$ by 
\[
\mho(H,M)=\left\{\left.(x_U)_{U\leq H}\in\prod_{U\leq H}\bz M(U)\,\right|\,
\con{U}{h}(x_U)=x_{\r{h}U}\quad\mbox{for all}\quad h\in H\right\}. 
\]
There is a ring homomorphism $\rho_H:\mathrm{\Omega}(H,M)\to\mho(H,M)$ given by 
\[
[(H/K)_s]\mapsto\left(\sum_{hK\in H/K,\,U\leq\r{h}\!K}\res{U}{\r{h}\!K}
(\r{h}\!s)\right)_{U\leq H} 
\]
for all $(K,s)\in\R(H,M)$ (cf. \cite[2.3]{Bol98}), which is called the mark 
homomorphism. 
\begin{definition}
\label{def:ghost}
We define a Green functor 
$\mho^M=(\mho^M,\con{}{},\res{}{},\ind{}{})$ for $G$ by 
\[
\begin{array}{l}
\vspace{0.3cm}\displaystyle
\mho^M(H)=\mho(H,M),\\
\vspace{0.3cm}
\con{H}{g}((x_U)_{U\leq H})=(\con{U}{g}(x_U))_{\r{g}\!U\leq\r{g}\!H}, \\
\vspace{0.3cm}
\res{K}{H}((x_U)_{U\leq H})=(x_U)_{U\leq K},\\
\displaystyle\ind{K}{H}((y_U)_{U\leq K})
=\sum_{hK\in H/K}(a_L^h)_{L\leq H} 
\end{array}
\]
for all $g\in G$, $K\leq H\leq G$, $(x_U)_{U\leq H}\in\mho(H,M)$, 
and $(y_U)_{U\leq K}\in\mho(K,M)$, where
\[
a_L^h=\left\{
\begin{array}{cl}
\vspace{0.1cm}
\con{U}{h}(y_U)\quad&\mbox{if }L=\r{h}U\mbox{ with }U\leq K,\\
0\quad&\mbox{otherwise}. 
\end{array}
\right.
\]
\end{definition}
\begin{remark}
\label{rem:gohst}
Let us remind the plus constructions $-_+:\ARes(G)\to\Green(G)$ and 
$-^+:\ACon(G)\to\Green(G)$, which are defined in \cite{Bol98}. The Green 
functors $\bz M_+$ and ${\bz M}^+$ are isomorphic to $\mathrm{\Omega}^M$ and 
$\mho^M$, respectively (cf. \cite[Proposition 3.1]{TY14}), and there is 
a morphism of Green functors $\rho:\mathrm{\Omega}^M\to\mho^M$ defined to be 
a family of the mark homomorphisms $\rho_H$ for $H\leq G$, which is called 
the mark morphism.
\end{remark}
\begin{definition}
\label{def:gho}
For each $H\leq G$, we define an additive group 
$\widetilde{\mathrm{\Omega}}(H,M)$ to be 
\[
\widetilde{\mathrm{\Omega}}(H,M)=\coprod_{(K,s)\in\R(H,M)}\bz, 
\]
and define an additive map 
$\varphi_H:\mathrm{\Omega}(H,M)\to\widetilde{\mathrm{\Omega}}(H,M)$ by 
\[
[(H/K)_s]\mapsto(|\{hK\in H/K\st U\leq\r{h}\!K\quad\mbox{and}\quad 
t=\res{U}{\r{h}\!K}(\r{h}\!s)\}|)_{(U,t)\in\R(H,M)}
\]
for all $(K,s)\in\R(H,M)$ (cf. \cite[p. 130]{TY14}). 
\end{definition}
\par
Let $H\leq G$. Given $(U,t)\in\SS(H,M)$ and 
$(x_{(K,s)})_{(K,s)\in\R(H,M)}\in\widetilde{\mathrm{\Omega}}(H,M)$, we set 
$x_{(U,t)}=x_{(K,s)}$ if $h.(U,t)=(K,s)\in\R(H,M)$ for some $h\in H$. 
\par
There exists an additive map 
$\varsigma_H:\widetilde{\mathrm{\Omega}}(H,M)\to\mathrm{\Omega}(H,M)$ given by 
\[
(x_{(K,s)})_{(K,s)\in\R(H,M)}\mapsto
\sum_{(U,t)\in\SS(H,M)}x_{(U,t)}\sum_{L\leq U}|L|
\mu(L,U)\left[(H/L)_{\res{L}{U}(t)}\right]
\]
for all $(x_{(K,s)})_{(K,s)\in\R(H,M)}\in\widetilde{\mathrm{\Omega}}(H)$, 
where $\mu$ is the M\"{o}bius function of the partially ordered set consisting 
of all subgroups of $G$ with the binary relation $\leq$ (cf. \cite{Aig79}). 
\begin{proposition}
\label{pro:exp}
For any $H\leq G$, 
\[
\varsigma_H\circ\varphi_H=|H|\cdot\id_{\mathrm{\Omega}(H,M)}
\qquad\mbox{and}\qquad 
\varphi_H\circ\varsigma_H=|H|\cdot\id_{\widetilde{\mathrm{\Omega}}(H,M)}. 
\]
\end{proposition}
{\itshape Proof.}
Let $H\leq G$. For each $(K,s)\in\R(H,M)$, we have 
\[
\begin{array}{l}
\vspace{0.1cm}
\varsigma_H\circ\varphi_H([(H/K)_s])\\
\vspace{0.1cm}\hspace{1.5cm}
\displaystyle=
\sum_{(U,t)\in\SS(H,M)}\varphi_{(U,t)}([(H/K)_s])
\sum_{L\leq U}|L|\mu(L,U)\left[(H/L)_{\res{L}{U}(t)}\right]\\
\vspace{0.1cm}\hspace{1.5cm}
\displaystyle=\sum_{hK\in H/K}\sum_{L\leq\r{h}\!K}
\sum_{L\leq U\leq\r{h}\!K}|L|\mu(L,U)
\left[(H/L)_{\res{L}{\r{h}\!K}(\r{h}\!s)}\right]\\
\hspace{1.5cm}
=|H|\cdot[(H/K)_s], 
\end{array}
\]
where 
\[
\varphi_{(U,t)}([(H/K)_s])=|\{hK\in H/K\st U\leq\r{h}\!K\quad\mbox{and}\quad
t=\res{U}{\r{h}\!K}(\r{h}\!s)\}|. 
\]
Thus $\varsigma_H\circ\varphi_H(x)=|H|\cdot x$ for any 
$x\in\mathrm{\Omega}(H,M)$. Let $\delta$ denote the Kronecker delta. For each 
$(x_{(K,s)})_{(K,s)\in\R(H,M)}\in\widetilde{\mathrm{\Omega}}(H,M)$, we have 
\[
\begin{array}{l}
\vspace{0.2cm}
\varphi_H\circ\varsigma_H((x_{(K,s)})_{(K,s)\in\R(H,M)})\\
\hspace{1.2cm}\displaystyle=
\sum_{(U,t)\in\SS(H,M)}x_{(U,t)}\sum_{L\leq U}|L|\mu(L,U)
\left(\sum_{h\in H/L,\,K\leq\r{h}\!L}
\delta_{s\,\res{K}{\r{h}U}(\r{h}t)}\right)_{(K,s)\in\R(H,M)}\\
\vspace{0.2cm}\hspace{1.2cm}\displaystyle=\sum_{h\in H}\sum_{U\leq H}
\left(\sum_{\r{h^{-1}}\!K\leq L\leq U}\mu(\r{h}\!L,\r{h}U)\sum_{t\in M(U)}
x_{(\r{h}U,\r{h}t)}\delta_{s\,\res{K}{\r{h}U}(\r{h}t)}\right)_
{(K,s)\in\R(H,M)}\\ 
\vspace{0.2cm}\hspace{1.2cm}\displaystyle=\sum_{h\in H}\sum_{U\leq H}
\left(\delta_{K\r{h}U}x_{(K,s)}\right)_{(K,s)\in\R(H,M)}\\
\hspace{1.2cm}\displaystyle
=|H|\cdot\left(x_{(K,s)}\right)_{(K,s)\in\R(H,M)}. 
\end{array}
\]
Hence $\varphi_H\circ\varsigma_H(y)=|H|\cdot y$ for any 
$y\in\widetilde{\mathrm{\Omega}}(H,M)$. This completes the proof. 
$\Box$
\par\bigskip
Given $H\leq G$ and $(K,s)\in\R(H,M)$, we set 
\[
N_H(K,s)=\{h\in H\st\r{h}\!K=K\quad\mbox{and}\quad\r{h}\!s=s\}, 
\]
which is a subgroup of the normalizer $N_H(K)$ of $K$ in $H$, and set 
\[
W_H(K,s)=N_H(K,s)/K. 
\]
\begin{proposition} 
\label{pro:kappa}
Let $H\leq G$. Given $(K,s)\in\R(H,M)$ and $L\leq H$, define 
\[
y^{(K,s)}_L:=\left\{
\begin{array}{cl}
\vspace{0.1cm}\displaystyle
\sum_{hN_H(K,s)\in N_H(K)/N_H(K,s)}\r{rh}\!s\quad
&\mbox{if }L=\r{r}\!K\mbox{ for some }r\in H,\\
0\quad&\mbox{otherwise}. 
\end{array}\right.
\]
The elements $(y^{(K,s)}_L)_{L\leq H}$ of $\mho(H,M)$ for 
$(K,s)\in\R(H,M)$ form a free $\bz$-basis of $\mho(H,M)$, and the additive map 
$\kappa_H:\widetilde{\mathrm{\Omega}}(H,M)\to\mho(H,M)$ given by 
\[
(\delta_{(K,s)\,(U,t)})_{(U,t)\in\R(H,M)}\mapsto(y^{(K,s)}_L)_{L\leq H}
\]
for all $(K,s)\in\R(H,M)$ is an isomorphism. 
Moreover, the diagram 
\begin{center}
\begin{picture}(90,90)(0,0)
\put(-8,60){\makebox(20,20)[c]{$\mathrm{\Omega}(H,M)$}}
\put(68,60){\makebox(20,20)[c]{$\widetilde{\mathrm{\Omega}}(H,M)$}}
\put(30,70){\makebox(20,20)[c]{\footnotesize{$\varphi_H$}}}
\put(30,70){\vector(1,0){20}}
\put(20,60){\vector(1,-1){40}}
\put(78,60){\vector(0,-1){40}}
\put(20,30){\makebox(20,20)[l]{\footnotesize{$\rho_H$}}}
\put(78,30){\makebox(20,20)[r]{\footnotesize{$\kappa_H$}}}
\put(68,0){\makebox(20,20)[c]{$\mho(H,M)$}}
\end{picture}
\end{center}
is commutative, and $\rho_H$ and $\varphi_H$ are injective. 
\end{proposition}
{\itshape Proof.}
The proof of the first assertion is straightforward (cf. \cite[\S9]{TY14}). We 
have 
\[
\begin{array}{l}
\vspace{0.3cm}
\kappa_H\circ\varphi_H([(H/K)_s])=\\
\vspace{0.1cm}\hspace{2cm}\displaystyle
\left(\sum_{(U,t)\in\R(H,M)}\sum_{hK\in H/K,\,U\leq\r{h}\!K}
\delta_{t\,\res{U}{\r{h}\!K}(\r{h}\!s)}y^{(U,t)}_L\right)_{L\leq H}\\
\vspace{0.1cm}\hspace{2cm}\displaystyle
\displaystyle=\left(\sum_{hK\in H/K}\sum_{(U,t)\in\R(H,M)}
\sum_{h_1N_H(U,t)\in N_H(U)/N_H(U,t)}z_L^{\{\r{h}\!s,\r{h_1}t\}}
\right)_{L\leq H}\\
\vspace{0.1cm}\hspace{2cm}\displaystyle
\displaystyle=\left(\sum_{hK\in H/K,\,L\leq\r{h}\!K}\res{L}{\r{h}\!K}
(\r{h}\!s)\right)_{L\leq H}\\
\hspace{2cm}\displaystyle=\rho_H([(H/K)_s]), 
\end{array}
\]
where 
\[
z^{\{\r{h}\!s,\r{h_1}t\}}_L=\left\{
\begin{array}{cl}
\vspace{0.1cm}\displaystyle
\r{rh_1}t\quad&\mbox{if }L=\r{r}U\leq\r{rh}\!K\mbox{ for some }r\in H
\mbox{ and if }t=\res{U}{\r{h}\!K}(\r{h}\!s),\\
0\quad&\mbox{otherwise}, 
\end{array}\right.
\]
for all $(K,s)\in\R(H,M)$, and thus $\kappa_H\circ\varphi_H=\rho_H$. Moreover, 
by Proposition \ref{pro:exp}, $\varphi_H$ is injective, and so is $\rho_H$ 
(cf. \cite[2.4]{Bol98}). This completes the proof.
$\Box$
\begin{proposition}
\label{pro:ghobasis}
For each $H\leq G$, 
\[
\widetilde{\mathrm{\Omega}}(H,M)
=\bigoplus_{(U,t)\in\R(H,M)}\frac{1}{\,|W_H(U,t)|\,}
\varphi_H([(H/U)_t])\bz. 
\]
\end{proposition}
{\itshape Proof.}
The proof is analogous to that of \cite[(80.15) Proposition]{CR87}. 
$\Box$
\begin{definition}
\label{def:obs}
For each $H\leq G$, we define an additive group $\obs(H,M)$ to be 
\[
\coprod_{(U,t)\in\R(H,M)}\bz/|W_H(U,t)|\bz. 
\]
\end{definition}
\begin{lemma}
\label{lem:gohstbasis}
For each $H\leq G$, 
\[
\widetilde{\mathrm{\Omega}}(H,M)/\im\varphi_H\cong\obs(H,M). 
\]
\end{lemma}
{\itshape Proof.}
By Propositions \ref{pro:monoid} and \ref{pro:kappa}, we have 
\[
\im\varphi_H=\bigoplus_{(U,t)\in\R(H,M)}\varphi_H([(H/U)_t])\bz. 
\]
Hence the assertion follows from Proposition \ref{pro:ghobasis}. 
$\Box$
\par\bigskip
Let $p$ be a prime, and let $\infty$ be just a symbol. For each $\bz$-module 
$R$, we set $R_{(p)}=\bz_{(p)}\otimes_\bz R$, where $\bz_{(p)}$ is 
the localization of $\bz$ at $p$, and $R_{(\infty)}=R$. 
\par
Let $H\leq G$. Given $(U,t)\in\SS(H,M)$, we denote by $W_H(U,t)_p$ a Sylow 
$p$-subgroup of $W_H(U,t)$, and set $W_H(U,t)_\infty=W_H(U,t)$. By 
Proposition \ref{pro:monoid}, the elements $[(G/K)_s]$ for $(K,s)\in\R(H,M)$ 
are supposed to form a free $\bz_{(p)}$-basis of $\mathrm{\Omega}(H,M)_{(p)}$. 
We identify $\widetilde{\mathrm{\Omega}}(H,M)_{(p)}$ and $\obs(H,M)_{(p)}$ 
with 
\[
\coprod_{(K,s)\in\R(H,M)}\bz_{(p)}\quad\mbox{and}\quad
\coprod_{(U,t)\in\R(H,M)}\bz_{(p)}/|W_G(U,t)_p|\bz_{(p)}, 
\]
respectively. Let $\varphi_H^{(p)}$ denote the $\bz_{(p)}$-module homomorphism 
from $\mathrm{\Omega}(H,M)_{(p)}$ to $\widetilde{\mathrm{\Omega}}(H,M)_{(p)}$ 
determined by $\varphi_H$. Then it follows from Lemma \ref{lem:gohstbasis} 
that 
\begin{equation}
\label{eq:ghostbasis}
\widetilde{\mathrm{\Omega}}(H,M)_{(p)}/\im\varphi_H^{(p)}\simeq
\obs(H,M)_{(p)}. 
\end{equation}
\par
Henceforth, we denote by $p$ a prime or the symbol $\infty$, and set 
$\varphi_H^{(\infty)}=\varphi_H$. The expression `$a\bmod b$' with 
$a,\,b\in\bz_{(p)}$ denotes the coset of $a+b\bz_{(p)}$ of $b\bz_{(p)}$ in 
$\bz_{(p)}$ containing $a$. Given $(U,t)\in\SS(H,M)$ and 
$(x_{(K,s)})_{(K,s)\in\R(H,M)}\in\widetilde{\mathrm{\Omega}}(H,M)_{(p)}$, we 
set $x_{(U,t)}=x_{(K,s)}$ if $h.(U,t)=(K,s)\in\R(H,M)$ for some $h\in H$. 
\par
We define a $\bz_{(p)}$-module homomorphism 
$\psi_H^{(p)}:\widetilde{\mathrm{\Omega}}(H,M)_{(p)}\to\obs(H,M)_{(p)}$ by 
\[
(x_{(K,s)})_{(K,s)\in\R(H,M)}\mapsto
\left(\sum_{\genfrac{}{}{0pt}{3}{rU\in W_H(U,t)_p,}{s\in M(\lr{r}U)_t}}
x_{(\lr{r}U,s)}\bmod|W_H(U,t)_p|\right)_{(U,t)\in\R(H,M)}, 
\]
where 
\[
M(\lr{r}U)_t=\{s\in M(\lr{r}U)\st\res{U}{\lr{r}U}(s)=t\}, 
\]
for all $(x_{(K,s)})_{(K,s)\in\R(H,M)}\in\widetilde{\mathrm{\Omega}}(H,M)$ 
(cf. \cite[p. 133]{TY14}). 
\begin{lemma}
\label{lem:CF}
For each $H\leq G$, $\psi_H^{(p)}$ is surjective. 
\end{lemma}
{\itshape Proof.}
The lemma follows from \cite[Lemma 9.3]{TY14}. (See also the proof of 
Lemma \ref{lem:psi}.) 
$\Box$
\par\bigskip
Let $H\leq G$, and let $(K,s)\in\R(H,M)$. We have  
\[
\begin{array}{l}
\vspace{0.1cm}\displaystyle
\psi_H^{(p)}\left(\varphi_H^{(p)}([(H/K)_s])\right)\\
\hspace{2.0cm}\displaystyle
=\left(\sum_{rU\in W_H(U,t)_p}
|\inv_{\lr{r}U}((H/K)_s)_{(U,t)}|\bmod|W_H(U,t)_p|\right)_{(U,t)\in\R(H,M)}, 
\end{array}
\]
where 
\[
\inv_{\lr{r}U}((H/K)_s)_{(U,t)}=\{hK\in H/K\st\lr{r}U\leq\r{h}\!K
\quad\mbox{and}\quad\res{U}{\r{h}\!K}(\r{h}\!s)=t\}. 
\]
Given $(U,t)\in\R(H,M)$ and $F\leq W_H(U,t)$, it follows from 
\cite[Lemma 9.2]{TY14} that 
\[
\sum_{rU\in F}|\inv_{\lr{r}U}((H/K)_s)_{(U,t)}|\equiv0\pmod{|F|}. 
\]
(See also the proof of Lemma \ref{lem:fundamental}). Hence we have 
\begin{equation}
\label{eq:Cauchy-Frobenius}
\im\varphi_H^{(p)}\subseteq\ker\psi_H^{(p)}. 
\end{equation}
\par
The following theorem is a generalization of \cite[(4.4) Theorem]{OY01}. 
\begin{theorem}[Fundamental theorem]
\label{thm:fnd}
For each $H\leq G$, the sequence 
\[
0\longrightarrow\mathrm{\Omega}(H,M)_{(p)}
\overset{\varphi_H^{(p)}}{\longrightarrow}
\widetilde{\mathrm{\Omega}}(H,M)_{(p)}\overset{\psi_H^{(p)}}{\longrightarrow}
\obs(H,M)_{(p)}\longrightarrow0
\]
of $\bz_{(p)}$-modules is exact. 
\end{theorem}
{\itshape Proof.}
The theorem is a special case of \cite[Theorem 9.4]{TY14}. We give a standard 
and concise proof of the theorem under the assumption that all the monoids 
$M(H)$ for $H\leq G$ are finite. Let $H\leq G$. By Proposition \ref{pro:kappa} 
and Lemma \ref{lem:CF}, it is enough to verify that 
$\im\varphi_H^{(p)}=\ker\psi_H^{(p)}$. Since $\psi_H^{(p)}$ is surjective, 
it follows that 
\[
\widetilde{\mathrm{\Omega}}(H,M)_{(p)}/\ker\psi_H^{(p)}\simeq\obs(H,M)_{(p)}. 
\]
Moreover, by Eqs.(\ref{eq:ghostbasis}) and (\ref{eq:Cauchy-Frobenius}), 
there is a sequence 
\[
\obs(H,M)_{(p)}\stackrel{\sim}{\to}
\widetilde{\mathrm{\Omega}}(H,M)_{(p)}/\im\varphi_H^{(p)}\to
\widetilde{\mathrm{\Omega}}(H,M)_{(p)}/\ker\psi_H^{(p)}  
\]
of finite groups, where the first arrow is an isomorphism and the second one 
is a natural surjection. Hence we have $\im\varphi_H^{(p)}=\ker\psi_H^{(p)}$, 
completing the proof. 
$\Box$
%----------------------------------------------------------------
%
%newsection
%
%----------------------------------------------------------------
\section{Lattice Burnside ring functors}
%----------------------------------------------------------------
Let $\L$ be a complete lattice with the binary relation $\leq$ (cf. 
\cite{Bir67}). For each subset $\Sigma$ of $\L$, $\inf\Sigma$ denotes 
the greatest lower bound of $\Sigma$ in $\L$, and $\sup\Sigma$ denotes 
the least upper bound of $\Sigma$ in $\L$. Given $x,\,y\in\L$, we set 
$x\wedge y=\inf\{x,\,y\}$ and $x\vee y=\sup\{x,\,y\}$. We consider $\L$ as 
a monoid with the binary operation given by 
\[
x\cdot y=x\wedge y
\]
for all $x,\,y\in\L$. (Of course, the binary operation $\wedge$ is 
associative.) The identity of the monoid $\L$ is the greatest element of $\L$. 
We call $\L$ a finite $G$-lattice if $\L$ is a finite left $G$-set and 
the binary relation $\leq$ is invariant under the action of $G$. Assume that 
$\L$ is a finite $G$-lattice. Then it turns out that $\L$ is a finite 
$G$-monoid (see Example \ref{epl:crossed}). Given $g\in G$ and $s\in \L$, 
$\r{g}\!s$ denotes the effect of $g$ on $s$. 
\begin{example}
\label{epl:subset}
The set $\Lambda$ of subsets of a finite left $G$-set is considered as 
a finite $G$-lattice with the binary relation given by inclusion and 
the action of $G$ given by 
\[
\r{g}\!I=\{gx\st x\in I\}
\]
for all $g\in G$ and $I\in\Lambda$; the binary operations $\wedge$ and 
$\vee$ are $\cap$ and $\cup$, respectively. 
\end{example}
\begin{proposition}
\label{pro:lattice}
Let $\L$ be a finite $G$-lattice. For each $H\leq G$, let $\L_H$ be 
a nonempty sublattice of $\L$. Suppose that the family of nonempty sublattices 
$\L_H$ of $\L$ for $H\leq G$ satisfies the following conditions. 
\def\theenumi{\arabic{enumi}}
\begin{enumerate}
\item
$\L_{\r{g}\!H}=\{\r{g}\!s\st s\in\L_H\}$ for any $g\in G$. 
\item
$\r{h}\!s=s$ for any $h\in H$ and $s\in\L_H$. 
\item
$s\wedge\sup\L_K\in\L_K$ for any $K\leq H$ and $s\in\L_H$. 
\item
$\sup\L_K\leq\sup\L_H$ for any $K\leq H$. 
\end{enumerate}
\def\theenumi{\alph{enumi}}
Then there exists a monoid functor $M_\L=(M_\L,\con{}{},\res{}{})$ for $G$ 
given by 
\[
\begin{array}{l}
\vspace{0.3cm}
M_\L(H)=\L_H, \\
\vspace{0.3cm}
\con{H}{g}:\L_H\to\L_{\r{g}\!H},\;s\mapsto\r{g}\!s,\\
\res{K}{H}:\L_H\to\L_K,\;s\mapsto s\wedge\sup\L_K
\end{array}
\]
for all $g\in G$ and $K\leq H\leq G$. 
\end{proposition}
{\itshape Proof.}
The proof is straightforward. 
$\Box$
\par\bigskip
Let $\L$ be a finite $G$-lattice, and let $M_\L=(M_\L,\con{}{},\res{}{})$ be 
the monoid functor given in Proposition \ref{pro:lattice}. The $M_\L$-Burnside 
ring functor is called the lattice Burnside ring functor on $M_\L$, and 
$\mathrm{\Omega}(G,M_\L)$ is called the lattice Burnside ring of $G$ 
associated to the family of nonempty sublattices $\L_H$ of $\L$ for $H\leq G$. 
\begin{example}
\label{epl:crossed-lattice}
If $\sup\L$ is $G$-invariant, then the $C_\L$-Burnside ring functor given in 
Example \ref{epl:crossed} with $S=\L$ is the lattice Burnside ring functor on 
$C_\L$. In particular, the $C_\Lambda$-Burnside ring functor, where $\Lambda$ 
is the finite $G$-lattice given in Example \ref{epl:subset}, is the lattice 
Burnside ring functor on $C_\Lambda$, because $\sup\Lambda$ is $G$-invariant. 
\end{example}
\par
Let $\SSS(G)$ denote the set of subgroups of $G$. We consider $\SSS(G)$ to be 
a finite $G$-lattice with the binary relation given by inclusion and 
the action of $G$ given by conjugation, and write $\SSS=\SSS(G)$ for the sake 
of shortness. 
\begin{example}
\label{epl:slice}
For each  $H\leq G$, we define a nonempty sublattice $\SSS_{\geq H}$ of $\SSS$ 
to be the set consisting of all subgroups of $G$ containing $H$. By 
Proposition \ref{pro:lattice}, there exists a monoid functor 
$M_\SSS=(M_\SSS,\con{}{},\res{}{})$ for $G$ given by 
\[
\begin{array}{l}
\vspace{0.3cm}
M_\SSS(H)=\SSS_{\geq H}, \\
\vspace{0.3cm}
\con{H}{g}:\SSS_{\geq H}\to\SSS_{\geq\r{g}\!H},\;E\mapsto\r{g}\!E,\\
\res{K}{H}:\SSS_{\geq H}\to\SSS_{\geq K},\;E\mapsto E
\end{array}
\]
for all $g\in G$ and $K\leq H\leq G$. Let $H\leq G$. Given 
$E,\,F\in\SSS_{\geq H}$, 
\[
E\cdot F=E\cap F. 
\]
By Proposition \ref{pro:monoid}, multiplication on 
$\mathrm{\Omega}(H,M_\SSS)$ is given by 
\[
[(H/K)_E]\cdot[(H/U)_F]
=\sum_{KhU\in K\backslash H/U}[(H/(K\cap\r{h}U))_{E\cap\r{h}F}]
\]
for all $(K,E),\,(U,F)\in\R(H,M_\SSS)$. The lattice Burnside ring 
$\mathrm{\Omega}(G,M_\SSS)$ is isomorphic to the slice Burnside ring 
$\mathrm{\Xi}(G)$ of $G$ introduced by S. Bouc \cite{Bou12}. 
\end{example}
\par
Almost all results on the ring structure of lattice Burnside rings are 
relative to the extension of 
the ring structure of slice Burnside rings. There is another example. 
\begin{example}
\label{epl:coslice}
For each $H\leq G$, we define a nonempty sublattice $\SSS_{\unlhd H}$ of 
$\SSS$ to be the set consisting of all normal subgroups of $H$. By 
Proposition \ref{pro:lattice}, there exists a monoid functor 
$M_\SSS^\circ=(M_\SSS^\circ,\con{}{},\res{}{})$ for $G$ given by 
\[
\begin{array}{l}
\vspace{0.3cm}
M_\SSS^\circ(H)=\SSS_{\unlhd H}, \\
\vspace{0.3cm}
\con{H}{g}:\SSS_{\unlhd H}\to\SSS_{\unlhd\r{g}\!H},\;L\mapsto\r{g}\!L,\\
\res{K}{H}:\SSS_{\unlhd H}\to\SSS_{\unlhd K},\;L\mapsto L\cap K
\end{array}
\]
for all $g\in G$ and $K\leq H\leq G$. Let $H\leq G$. Given 
$L,\,N\in\SSS_{\unlhd H}$, 
\[
L\cdot N=L\cap N. 
\]
By Proposition \ref{pro:monoid}, multiplication on 
$\mathrm{\Omega}(H,M_\SSS^\circ)$ is given by 
\begin{equation}
\label{eq:product}
[(H/K)_L]\cdot[(H/U)_N]
=\sum_{KhU\in K\backslash H/U}[(H/(K\cap\r{h}U))_{L\cap\r{h}N}]
\end{equation}
for all $(K,L),\,(U,N)\in\R(H,M_\SSS^\circ)$. 
\end{example}
\par
The lattice Burnside ring $\mathrm{\Omega}(G,M_\L)$ is isomorphic to 
an abstract Burnside ring (see Theorem \ref{thm:trans}). But, the fundamental 
theorem of $\mathrm{\Omega}(G,M_\L)$ (see Theorem \ref{thm:fnd}) is not 
derived from that of an abstract Burnside ring (see 
Theorem \ref{thm:abstract}). So we need to make an adjustment of 
Theorem \ref{thm:fnd} with $M=M_\L$. Recall that $p$ is a prime or the symbol 
$\infty$. For each $H\leq G$, we consider 
$\widetilde{\mathrm{\Omega}}(H,M_\L)_{(p)}$ as the ring 
\[
\prod_{(K,s)\in\R(H,M)}\bz_{(p)}, 
\]
and provide complementary maps from 
$\widetilde{\mathrm{\Omega}}(H,M_\L)_{(p)}$ to itself. 
\begin{definition}
\label{def:ring}
Given $U\leq G$, let $\mu_U$ denote the M\"{o}bius function of the partially 
ordered set $\L_U$ with the binary relation $\leq$. For each $H\leq G$, we 
define maps $\alpha_H^{(p)}:\widetilde{\mathrm{\Omega}}(H,M_\L)_{(p)}\to
\widetilde{\mathrm{\Omega}}(H,M_\L)_{(p)}$ and 
$\beta_H^{(p)}:\widetilde{\mathrm{\Omega}}(H,M_\L)_{(p)}\to
\widetilde{\mathrm{\Omega}}(H,M_\L)_{(p)}$ by 
\[
\alpha_H^{(p)}((x_{(K,s)})_{(K,s)\in\R(H,M_\L)})
=\left(\sum_{t\leq s\in\L_U}x_{(U,s)}\right)_{(U,t)\in\R(H,M_\L)}
\]
and
\[
\beta_H^{(p)}((x_{(K,s)})_{(K,s)\in\R(H,M_\L)})=\left(
\sum_{t\leq s\in\L_U}\mu_U(t,s)x_{(U,s)}\right)_{(U,t)\in\R(H,M_\L)} 
\]
for all $(x_{(K,s)})_{(K,s)\in\R(H,M_\L)}\in
\widetilde{\mathrm{\Omega}}(H,M_\L)_{(p)}$, respectively. 
\end{definition}
\begin{lemma}
\label{lem:lattice}
For each $H\leq G$, $\beta_H^{(p)}\circ\alpha_H^{(p)}
=\alpha_H^{(p)}\circ\beta_H^{(p)}
=\id_{\widetilde{\mathrm{\Omega}}(H,M_\L)_{(p)}}$. 
\end{lemma}
{\itshape Proof.}
Let $H\leq G$, and let $(x_{(K,s)})_{(K,s)\in\R(H,M_\L)}\in
\widetilde{\mathrm{\Omega}}(H,M_\L)_{(p)}$. We have 
\[
\begin{array}{l}
\vspace{0.1cm}
\beta_H^{(p)}\circ\alpha_H^{(p)}((x_{(K,s)})_{(K,s)\in\R(H,M_\L)})\\
\vspace{0.1cm}\hspace{2cm}\displaystyle
=\left(\sum_{t\leq s_2\in\L_U}\mu_U(t,s_2)\sum_{s_2\leq s_1\in\L_U}x_{(U,s_1)}
\right)_{(U,t)\in\R(H,M_\L)}\\
\vspace{0.1cm}
\hspace{2cm}\displaystyle
=\left(\sum_{t\leq s_1\in\L_U}x_{(U,s_1)}
\sum_{t\leq s_2\in\L_U^{\leq s_1}}\mu_H(t,s_2)\right)_{(U,t)\in\R(H,M_\L)}\\
\hspace{2cm}
=(x_{(U,t)})_{(U,t)\in\R(H,M_\L)}, 
\end{array}
\]
where $\L_U^{\leq s_1}=\{s_2\in\L_U\st s_2\leq s_1\}$, and 
\[
\begin{array}{l}
\vspace{0.1cm}
\alpha_H^{(p)}\circ\beta_H^{(p)}((x_{(K,s)})_{(K,s)\in\R(H,M_\L)})\\
\vspace{0.1cm}\hspace{2cm}\displaystyle
=\left(\sum_{t\leq s_2\in\L_U}\sum_{s_2\leq s_1\in\L_U}
\mu_U(s_2,s_1)x_{(U,s_1)}\right)_{(U,t)\in\R(H,M_\L)}\\
\vspace{0.1cm}
\hspace{2cm}\displaystyle
=\left(\sum_{t\leq s_1\in\L_U}x_{(U,s_1)}
\sum_{t\leq s_2\in\L_U^{\leq s_1}}
\mu_U(s_2,s_1)\right)_{(U,t)\in\R(H,M_\L)}\\
\hspace{2cm}
=(x_{(U,t)})_{(U,t)\in\R(H,M_\L)}. 
\end{array}
\]
This completes the proof. 
$\Box$
\begin{definition}
\label{def:psi}
Let $H\leq G$. Given $(U,t)\in\R(H,M)$ and $rU\in W_H(U,t)$, set 
\[
\L_{\lr{r}U}^{\geq t}=\{s\in\L_{\lr{r}U}\st s\geq t\}. 
\]
We define a $\bz_{(p)}$-module homomorphism 
$\widetilde{\psi}_H^{(p)}:\widetilde{\mathrm{\Omega}}(H,M_\L)_{(p)}
\to\obs(H,M_\L)_{(p)}$ by 
\[
(x_{(K,s)})_{(K,s)\in\R(H,M_\L)}\mapsto\left(
\sum_{\genfrac{}{}{0pt}{1}{rU\in W_H(U,t)_p,}{s\in\L_{\lr{r}U}^{\geq t}}}
x_{(\lr{r}U,s)}\bmod|W_H(U,t)_p|\right)_{(U,t)\in\R(H,M)}. 
\]
\end{definition}
\par
The proof of the following lemma is analogous to that of 
\cite[Lemma 4.3]{TY10}. 
\begin{lemma}
\label{lem:psi}
For each $H\leq G$, $\widetilde{\psi}_H^{(p)}$ is surjective. 
\end{lemma}
{\itshape Proof.}
Set $\ov{\delta}_{(K,s)}=(\delta_{(K,s)\;(U,t)}
\bmod|W(U,t)|_p)_{(U,t)\in\R(H,M_\L)}\in\obs(H,M_\L)_{(p)}$ for each 
$(K,s)\in\R(H,M_\L)$. Obviously, the elements 
$\ov{\delta}_{(K,s)}$ for $(K,s)\in\R(H,M_\L)$ form a $\bz_{(p)}$-basis of 
$\obs(H,M_\L)_{(p)}$. Now set 
\[
\R_0=\{(K,s)\in\R(H,M_\L)\st\ov{\delta}_{(K,s)}\not\in
\im\widetilde{\psi}_H^{(p)}\}. 
\]
We define a partially order $\leq_H$ on $\R(H,M_\L)$ by the rule that 
\[
(U,t)\leq_H(K,s)\quad\Longleftrightarrow\quad 
U\leq\r{h}\!K\quad\mbox{and}\quad t\leq\r{h}\!s\quad\mbox{for some}\quad 
h\in H. 
\]
Suppose that $\R_0\not=\emptyset$, and let $(K,s)$ be a minimal element of 
$\R_0$ with respect to $\leq_H$. Then no element $(U,t)$ of $\R_0-\{(K,s)\}$ 
satisfies $(U,t)\leq_H(K,s)$, and thus 
\[
\widetilde{\psi}_H^{(p)}((\delta_{(K,s)\,(U,t)})_{(U,t)\in\R(H,M)})
=(y_{(U,t)})_{(U,t)\in\R(H,M)}, 
\]
where
\[
y_{(U,t)}=\left\{
\begin{array}{ll}
\vspace{0.1cm}
1\bmod|W_H(U,t)_p|\quad\mbox{if}\quad(U,t)=(K,s),\\
0\bmod|W_H(U,t)_p|\quad\mbox{if}\quad(U,t)\in\R_0-\{(K,s)\}. 
\end{array}
\right.
\]
But, $\ov{\delta}_{(U,t)}\in\im\widetilde{\psi}_H^{(p)}$ for any 
$(U,t)\not\in\R_0$, which yields 
$\ov{\delta}_{(K,s)}\in\im\widetilde{\psi}_H^{(p)}$. 
This is a contradiction. Consequently, we have $\R_0=\emptyset$, completing 
the proof. 
$\Box$
\par\bigskip
The proof of the following lemma is analogous to that of 
\cite[Lemma 9.2]{TY14}. 
\begin{lemma}
\label{lem:fundamental}
Let $H\leq G$, and let $(K,s),\,(U,t)\in\R(H,M_\L)$. Set 
\[
\inv_U((H/K)_s)_{\geq t}
=\{hK\in H/K\st U\leq\r{h}\!K\quad\mbox{and}\quad t\leq\r{h}\!s\}. 
\]
Then for any $F\leq W_H(U,t)$, 
\[
\sum_{rU\in F}|\{hK\in
\inv_U((H/K)_s)_{\geq t}\st\lr{r}U\leq\r{h}\!K\}|\equiv0\pmod{|F|}. 
\]
\end{lemma}
{\itshape Proof.}
Let $F\leq W_H(U,t)$, and make $\inv_U((H/K)_s)_{\geq t}$ into a left $F$-set 
by defining 
\[
rUhK=rhK
\]
for all $rU\in F$ and $hK\in\inv_U((H/K)_s)_{\geq t}$. Then for any $rU\in F$, 
\[
\{hK\in\inv_U((H/K)_s)_{\geq t}\st\lr{r}U\leq\r{h}\!K\}
=\{hK\in\inv_U((H/K)_s)_{\geq t}\st rUhK=hK\}, 
\]
which is the set of $F$-invariants in $\inv_U((H/K)_s)_{\geq t}$. Let 
$F\backslash\inv_U((H/K)_s)_{\geq t}$ denote the set of $F$-orbits in 
$\inv_U((H/K)_s)_{\geq t}$, By \cite[Lemma 2.7]{Yo90b}, we have 
\[
\begin{array}{ll}
\vspace{0.1cm}
\displaystyle
\sum_{rU\in F}|\{hK\in\inv_U((H/K)_s)_{\geq t}\st\lr{r}U\leq\r{h}\!K\}|
&=|F|\cdot|F\backslash\inv_U((H/K)_s)_{\geq t}| \\ 
\mbox{}&\equiv0\pmod{|F|}, 
\end{array}
\]
completing the proof. 
$\Box$
\par\bigskip
Let $H\leq G$. Given $x\in\mathrm{\Omega}(H,M_\L)_{(p)}$, it follows from 
Lemmas \ref{lem:lattice} and \ref{lem:fundamental} that 
\begin{equation}
\label{eq:psi-vsrphi}
(\widetilde{\psi}_H^{(p)}\circ\beta_H^{(p)})\circ
(\alpha_H^{(p)}\circ\varphi_H^{(p)})(x)
=\widetilde{\psi}_H^{(p)}\circ\varphi_H^{(p)}(x)=0\bmod|W_H(U,t)_p|. 
\end{equation}
We define a $\bz_{(p)}$-module homomorphism 
$\widetilde{\alpha}_H^{(p)}:\mathrm{\Omega}(H,M_\L)_{(p)}\to
\widetilde{\mathrm{\Omega}}(H,M_\L)_{(p)}$ by 
\[
[(H/K)_s]\mapsto(|\{hK\in H/K\st U\leq\r{h}\!K\quad\mbox{and}\quad 
t\leq\r{h}\!s\}|)_{(U,t)\in\R(H,M_\L)}
\]
for all $(K,s)\in\R(H,M_\L)$. Given $(U,t)\in\R(H,M_\L)$ and 
$rU\in W_H(U,t)_p$, set 
\[
s_{(r,t)}=\inf\L_{\lr{r}U}^{\geq t}. 
\]
We define a $\bz_{(p)}$-module homomorphism $\widetilde{\beta}_H^{(p)}:
\widetilde{\mathrm{\Omega}}(H,M_\L)_{(p)}\to\obs(H,M_\L)_{(p)}$ by 
\[
(x_{(K,s)})_{(K,s)\in\R(H,M)}\mapsto\left(\sum_{rU\in W_H(U,t)_p}
x_{(\lr{r}U,s_{(r,t)})}\bmod|W_H(U,t)_p|\right)_{(U,t)\in\R(H,M_\L)}
\]
for all 
$(x_{(K,s)})_{(K,s)\in\R(H,M_\L)}\in\widetilde{\mathrm{\Omega}}(H,M_\L)_{(p)}$. \begin{theorem}
\label{thm:fundamental}
For each $H\leq G$, the sequence 
\begin{equation}
\label{eq:lattice}
0\longrightarrow\mathrm{\Omega}(H,M_\L)_{(p)}
\overset{\alpha_H^{(p)}\circ\varphi_H^{(p)}}
{\longrightarrow}\widetilde{\mathrm{\Omega}}(H,M_\L)_{(p)}
\overset{\widetilde{\psi}_H^{(p)}\circ\beta_H^{(p)}}{\longrightarrow}
\obs(H,M_\L)_{(p)}\longrightarrow0
\end{equation}
of $\bz_{(p)}$-modules is exact; moreover, 
$\widetilde{\alpha}_H^{(p)}=\alpha_H^{(p)}\circ\varphi_H^{(p)}$, and 
$\widetilde{\beta}_H^{(p)}=\widetilde{\psi}_H^{(p)}\circ\beta_H^{(p)}$. 
\end{theorem}
{\itshape Proof.}
Let $H\leq G$. By Theorem \ref{thm:fnd}, Lemma \ref{lem:psi}, and 
Eq.(\ref{eq:psi-vsrphi}), the sequence 
\[
0\longrightarrow\mathrm{\Omega}(H,M_\L)_{(p)}\overset{\varphi_H^{(p)}}
{\longrightarrow}\widetilde{\mathrm{\Omega}}(H,M_\L)_{(p)}
\overset{\widetilde{\psi}_H^{(p)}}{\longrightarrow}\obs(H,M_\L)_{(p)}
\longrightarrow0
\]
of $\bz_{(p)}$-modules is exact (see the proof of Theorem \ref{thm:fnd}). 
Thus it follows from Lemma \ref{lem:lattice} that Eq.(\ref{eq:lattice}) is 
exact. Obviously, 
$\widetilde{\alpha}_H^{(p)}=\alpha_H^{(p)}\circ\varphi_H^{(p)}$. Suppose that 
\[
\widetilde{\psi}_H^{(p)}\circ\beta_H^{(p)}((x_{(K,s)})_{(K,s)\in\R(H,M)})
=(y_{(U,t)})_{(U,t)\in\R(H,M)}\in\obs(H,M_\L)_{(p)} 
\]
with $y_{(U,t)}\in\bz_{(p)}/|W_H(U,t)_p|\bz_{(p)}$. Then for each 
$(U,t)\in\R(H,M)$, 
\[
\begin{array}{ll}
\vspace{0.1cm}
y_{(U,t)}&\displaystyle
=\sum_{\genfrac{}{}{0pt}{3}{rU\in W_H(U,t)_p,}{s_1\in\L_{\lr{r}U}^{\geq t}}}
\sum_{s_1\leq s_2\in\L_{\lr{r}U}}\mu_{\lr{r}U}(s_1,s_2)x_{(\lr{r}U,s_2)}
\bmod|W_H(U,t)_p|\\
\vspace{0.1cm}
\mbox{}&\displaystyle
=\sum_{\genfrac{}{}{0pt}{3}{rU\in W_H(U,t)_p,}{s_2\in\L_{\lr{r}U}}}
\sum_{\genfrac{}{}{0pt}{3}{s_1\in\L_{\lr{r}U}^{\geq t},}{s_1\leq s_2}}
\mu_{\lr{r}U}(s_1,s_2)x_{(\lr{r}U,s_2)}\bmod|W_H(U,t)_p|\\
\mbox{}&\displaystyle
=\sum_{rU\in W_H(U,t)_p}x_{(\lr{r}U,s_{(r,t)})}
\bmod|W_H(U,t)_p|. 
\end{array}
\]
Hence we have 
$\widetilde{\beta}_H^{(p)}=\widetilde{\psi}_H^{(p)}\circ\beta_H^{(p)}$, 
completing the proof. 
$\Box$
\par\bigskip
For each $H\leq G$, the primitive idempotents of 
$\bq\otimes_\bz\mathrm{\Omega}(H)$ are the elements 
\[
e_K^{(H)}:=\frac{1}{\,|N_H(K)|\,}\sum_{U\leq K}|U|\mu(U,K)[H/U]
\]
of $\bq\otimes_\bz\mathrm{\Omega}(H)$ for $K\in\C(H)$ which is given in 
\cite{Glu81,Yo83a}. In addition, the primitive idempotents of 
$\bq\otimes_\bz\mathrm{\Xi}(G)$ are given in \cite[Corollary 5.5]{Bou12}, 
where $\mathrm{\Xi}(G)$ is the slice Burnside ring of $G$ (see 
Example \ref{epl:slice}). We are now successful in finding the primitive 
idempotents of $\bq\otimes_\bz\mathrm{\Omega}(H,M_\L)$ for $H\leq G$. 
\begin{theorem}
\label{thm:ring homomorphism}
Let $H\leq G$, and define a map 
$\zeta_H:\mho(H,M_\L)\to\widetilde{\mathrm{\Omega}}(H,M_\L)$ 
by 
\[
\left(\sum_{s\in\L_K}\ell_{(K,s)}s\right)_{K\leq G}\mapsto
\left(\sum_{t\leq s\in\L_U}\ell_{(U,s)}\right)_{(U,t)\in\R(H,M_\L)}
\]
for all $\ell_{(K,s)}\in\bz$ with $(K,s)\in\SS(H,M_\L)$. Let $\alpha_H$ 
denote $\alpha_H^{(\infty)}$. Then the diagram 
\begin{center}
\begin{picture}(90,90)(20,0)
\put(-8,60){\makebox(20,20)[c]{$\mathrm{\Omega}(H,M_\L)$}}
\put(68,60){\makebox(20,20)[c]{$\widetilde{\mathrm{\Omega}}(H,M_\L)$}}
\put(144,60){\makebox(20,20)[c]{$\widetilde{\mathrm{\Omega}}(H,M_\L)$}}
\put(30,70){\makebox(20,20)[c]{\footnotesize{$\varphi_H$}}}
\put(106,70){\makebox(20,20)[c]{\footnotesize{$\alpha_H$}}}
\put(30,70){\vector(1,0){20}}
\put(106,70){\vector(1,0){20}}
\put(20,60){\vector(1,-1){40}}
\put(96,20){\vector(1,1){40}}
\put(20,30){\makebox(20,20)[l]{\footnotesize{$\rho_H$}}}
\put(120,30){\makebox(20,20)[r]{\footnotesize{$\zeta_H$}}}
\put(68,0){\makebox(20,20)[c]{$\mho(H,M_\L)$}} 
\end{picture}
\end{center}
is commutative. Moreover, the map 
$\alpha_H\circ\varphi_H:\mathrm{\Omega}(H,M_\L)\to
\widetilde{\mathrm{\Omega}}(H,M_\L)$ is a ring monomorphism, and 
the primitive idempotents of $\bq\otimes_\bz\mathrm{\Omega}(H,M_\L)$ 
are the elements 
\[
e_{(K,s)}^{(H)}:=\frac{1}{\,|N_H(K,s)|\,}
\sum_{t\in\L_K^{\leq s}}\mu_K(t,s)\sum_{U\leq K}|U|\mu(U,K)
[(H/U)_{t\wedge\sup\L_U}], 
\]
where $\L_K^{\leq s}=\{t\in\L_K\st t\leq s\}$, for $(K,s)\in\R(H,M_\L)$. 
\end{theorem}
{\itshape Proof.}
By Proposition \ref{pro:kappa}, $\kappa_H\circ\varphi_H=\rho_H$. Moreover, 
$\alpha_H=\zeta_H\circ\kappa_H$, because 
\[
\begin{array}{l}
\vspace{0.1cm}
\alpha_H((\delta_{(K,s)\,(U,t)})_{(U,t)\in\R(H,M)})\\
\vspace{0.1cm}\hspace{1cm}
=(\delta_{K\,U}|\{s_1\in\L_K\st t\leq s_1\mbox{ and }\r{h}\!s=s_1
\mbox{ for some }h\in N_H(K)\}|)_{(U,t)\in\R(H,M)}\\
\hspace{1cm}=\zeta_H\circ\kappa_H((\delta_{(K,s)\,(U,t)})_{(U,t)\in\R(H,M)}) 
\end{array}
\]
for all $(K,s)\in\R(H,M_\L)$. Hence 
$\alpha_H\circ\varphi_H=\zeta_H\circ(\kappa_H\circ\varphi_H)
=\zeta_H\circ\rho_H$. Obviously, $\zeta_H$ is a ring isomorphism. Since 
$\rho_H$ is a ring monomorphism (cf. \cite[2.4]{Bol98}), so is 
$\alpha_H\circ\varphi_H$. This, combined with Proposition \ref{pro:exp} and 
Lemma \ref{lem:lattice}, shows that the primitive idempotents of 
$\bq\otimes_\bz\mathrm{\Omega}(H,M_\L)$ are the elements 
\[
\frac{1}{\,|H|\,}\,
\varsigma_H\circ\beta_H((\delta_{(K,s)\,(U,t)})_{(U,t)\in\R(H,M_\L)}) 
\]
for $(K,s)\in\R(H,M_\L)$, where $\beta_H=\beta_H^{(\infty)}$. Given 
$(K,s)\in\R(H,M_\L)$, we have 
\[
\beta_H((\delta_{(K,s)\,(U,t)})_{(U,t)\in\R(H,M_\L)})
=\left(\delta_{K\,U}\sum_{t\leq s_1\in\L_K^{(s)}}
\mu_K(t,s_1)\right)_{(U,t)\in\R(H,M_\L)} 
\]
where $\L_K^{(s)}=\{s_1\in\L_K\st s_1=\r{h}\!s\mbox{ for some }h\in N_H(K)\}$, 
and 
\[
\begin{array}{l}
\vspace{0.3cm}
\varsigma_H\circ\beta_H((\delta_{(K,s)\,(U,t)})_{(U,t)\in\R(H,M_\L)})\\
\vspace{0.1cm}\hspace{1cm}\displaystyle
=\frac{|H|}{\,|N_H(K)|\,}\sum_{t\in\L_K^{\leq(s)}}
\sum_{t\leq s_1\in\L_K^{(s)}}\mu_K(t,s_1)\sum_{U\leq K}|U|
\mu(U,K)[(H/U)_{t\wedge\sup\L_U}]\\
\vspace{0.1cm}\hspace{1cm}\displaystyle
=\frac{|H|}{\,|N_H(K)|\,}\sum_{s_1\in\L_K^{(s)}}
\sum_{t\in\L_K^{\leq s_1}}\mu_K(t,s_1)\sum_{U\leq K}|U|
\mu(U,K)[(H/U)_{t\wedge\sup\L_U}]\\
\hspace{1cm}\displaystyle
=\frac{|H|}{\,|N_H(K,s)|\,}\sum_{t\in\L_K^{\leq s}}
\mu_K(t,s)\sum_{U\leq K}|U|\mu(U,K)[(H/U)_{t\wedge\sup\L_U}], 
\end{array}
\]
where 
$\L_K^{\leq(s)}=\{t\in\L_K\st t\leq\r{h}\!s\mbox{ for some }h\in N_H(K)\}$. 
This completes the proof. 
$\Box$
\par\bigskip
The exact sequence (\ref{eq:lattice}) with $H=G$ is derived from the exact 
sequence (\ref{eq:abstract}) for abstract Burnside rings (see the next 
section). Also, the primitive idempotents $e_{(K,s)}^{(G)}$ of 
$\bq\otimes_\bz\mathrm{\Omega}(H,M_\L)$ for $(K,s)\in\R(H,M_\L)$, which are 
given in Theorem \ref{thm:ring homomorphism}, are explained in terms of 
the M\"{o}bius function of $\SS(G,M_\L)$ (see Remark \ref{rem:Mebius}). 
%----------------------------------------------------------------
%
%newsection
%
%----------------------------------------------------------------
\section{Units of lattice Burnside rings }
%----------------------------------------------------------------
We show that any lattice Burnside ring is isomorphic to an abstract Burnside 
ring, and explore the units of a lattice Burnside ring. 
\par
Let $\boldGamma$ be an essentially finite category, that is, a category 
equivalent to a finite category. The set of isomorphism classes of objects of 
$\boldGamma$ is denoted by $\boldGamma\!/\!\simeq$. Given objects $I$ and $J$ 
of $\boldGamma$, the set of morphisms from $I$ to $J$ is denoted by 
$\hom(I,J)$ or $\boldGamma(I,J)$. Since $\boldGamma$ is essentially finite, it 
follows that $\boldGamma\!/\!\simeq$ and $\boldGamma(I,J)$ are finite sets. 
\par
As before, $p$ denotes a prime or the symbol $\infty$. Let 
$\bz\boldGamma$ be the free abelian group on $\boldGamma\!/\!\simeq$, and set 
$\bz_{(p)}\boldGamma=\bz_{(p)}\otimes_\bz\bz\boldGamma$. The isomorphism class 
containing an object $I$ of $\boldGamma$ is denoted by $[I]$. The ghost ring 
$\bz_{(p)}^\boldGamma$ of $\bz_{(p)}\boldGamma$ is defined to be 
\[
\bz_{(p)}^\boldGamma=\prod_{[I]\in\boldGamma\!/\!\simeq}\bz_{(p)}. 
\]
We define a $\bz_{(p)}$-module homomorphism 
$\varphi_\boldGamma^{(p)}:\bz_{(p)}\boldGamma\to\bz_{(p)}^\boldGamma$ by 
\begin{equation}
\label{eq:varphi}
[J]\mapsto(|\boldGamma(I,J)|)_{[I]\in\boldGamma\!/\!\simeq}
\end{equation}
for all $[J]\in\boldGamma\!/\!\simeq$, and call it the Burnside map. 
\par
The $\bz_{(p)}$-module $\bz_{(p)}\boldGamma$ is called an abstract Burnside 
ring if it has a $\bz_{(p)}$-algebra structure such that 
$\varphi_\boldGamma^{(p)}$ is an injective $\bz_{(p)}$-algebra homomorphism 
(cf. \cite{YOT18}). 
\par
Let $I$ be an object of $\boldGamma$. We denote by $\aut(I)$ the group of 
automorphisms of $I$. For each $\sigma\in\aut(I)$, a morphism 
$c_\sigma:I\to I/\sigma$ in $\boldGamma$ is said to be a coequalizer of 
$\sigma$ and $\id_I$ if $c_\sigma\circ\sigma=c_\sigma$ 
and it is universal with this property; that is, for any morphism $f:I\to J$ 
in $\boldGamma$ which satisfies that $f\circ\sigma=f$, there exists a unique 
morphism $f_1:I/\sigma\to J$ in $\boldGamma$ such that $f=f_1\circ c_\sigma$: 
\begin{center}
\begin{picture}(90,90)(40,0)
\put(12,60){\makebox(20,20)[c]{$I$}}
\put(68,60){\makebox(20,20)[c]{$I$}}
\put(124,60){\makebox(20,20)[c]{$I/\sigma$}}
\put(40,71){\makebox(20,20)[c]{\footnotesize{$\id_I$}}}
\put(40,50){\makebox(20,20)[c]{\footnotesize{$\sigma$}}}
\put(94,66){\makebox(20,20)[c]{\footnotesize{$c_\sigma$}}}
\put(32,73){\vector(1,0){35}}
\put(32,67){\vector(1,0){35}}
\put(86,70){\vector(1,0){35}}
\put(86,60){\vector(1,-1){35}}
\dashline{1}(134,60)(134,57)(134,54)(134,51)(134,48)(134,45)(134,42)(134,39)(134,36)(134,33)(134,30)
\put(134,27){\vector(0,-1){2}}
\put(72,34){\makebox(20,20)[r]{\footnotesize{$f$}}}
\put(128,34){\makebox(20,20)[r]{\footnotesize{$f_1$}}}
\put(124,6){\makebox(20,20)[c]{$J$.}} 
\end{picture}
\end{center}
\par
Given an object $I$ of $\boldGamma$, let $\aut(I)_p$ denote 
a Sylow $p$-subgroup of $\aut(I)$, and let $\aut(I)_\infty$ denote $\aut(I)$. 
For an epi-mono factorization system of the category $\boldGamma$, we refer to 
\cite[1.4]{YOT18}. The assertions of the following theorem are presented in 
\cite[Theorems 2.4 and 2.6]{YOT18} under a slightly weaker assumption. 
\begin{theorem}[Fundamental theorem]
\label{thm:abstract}
Assume that $\boldGamma$ has an epi-mono factorization system and that for 
any object $I$ of $\boldGamma$, any element $\sigma$ of $\aut(I)_p$ has 
a coequalizer $c_\sigma:I\to I/\sigma$ of $\sigma$ and $\id_I$. Then there 
exists an exact sequence 
\begin{equation}
\label{eq:abstract}
0\longrightarrow\bz_{(p)}\boldGamma
\overset{\varphi_\boldGamma^{(p)}}{\longrightarrow}
\bz_{(p)}^\boldGamma\overset{\psi_\boldGamma^{(p)}}{\longrightarrow}
\obs_{(p)}(\boldGamma)
\longrightarrow0 
\end{equation}
of $\bz_{(p)}$-modules, where $\obs_{(p)}(\boldGamma)$ is the group of 
obstruction of $\boldGamma$ defined to be 
\[
\obs_{(p)}(\boldGamma)=\prod_{[I]\in\boldGamma\!/\!\simeq}
\bz_{(p)}/|\aut(I)_p|\bz_{(p)}, 
\]
and $\psi_\boldGamma^{(p)}:\bz_{(p)}^\boldGamma\to\obs_{(p)}(\boldGamma)$ is 
the Cauchy-Frobenius map given by 
\[
(x_I)_{[I]\in\boldGamma\!/\!\simeq}\mapsto\left(
\sum_{\sigma\in\aut(I)_p}x_{I/\sigma}\bmod|\aut(I)_p|
\right)_{[I]\in\boldGamma\!/\!\simeq}
\]
for all $(x_I)_{[I]\in\boldGamma\!/\!\simeq}\in\bz_{(p)}^\boldGamma$. 
Moreover, $\bz_{(p)}\boldGamma$ is an abstract Burnside ring. 
\end{theorem}
\par
If $\boldGamma$ satisfies the assumptions of Theorem {\rm\ref{thm:abstract}}, 
we consider $\bz_{(p)}\boldGamma$ to be the abstract Burnside ring which has 
a $\bz_{(p)}$-algebra structure such that $\varphi_\boldGamma^{(p)}$ is 
an injective $\bz_{(p)}$-algebra homomorphism. The origin of abstract Burnside 
rings is $\mathrm{\Omega}(G)$. 
\begin{example}
\label{epl:trans}
Let $\trans^G_\SSS$ be the category of transitive $G$-sets $G/H$ for $H\leq G$ 
and $G$-equivariant maps. Given $U,\,H\leq G$, we have 
\[
\hom(G/U,G/H)=\{\lambda_g:G/U\to G/H,\,rU\mapsto rgH\st gH\in\inv_U(G/H)\}, 
\]
where $\inv_U(G/H)=\{gH\in G/H\st U\leq\r{g}\!H\}$, so that $G/U\simeq G/H$ if 
and only if $U$ is a conjugate of $H$ (cf. \cite[(80.5) Proposition]{CR87}). 
Let $H\leq G$. Then 
\[
\aut(G/H)=\{\sigma_g:G/H\to G/H,\,rH\mapsto rgH\st gH\in N_G(H)/H\}. 
\]
Given $g\in N_G(H)$, there exists a coequalizer 
$c_{\sigma_g}:G/H\to (G/H)/\sigma_g$ of $\sigma_g$ and $\id_{G/H}$ given by 
$(G/H)/\sigma_g=G/(\lr{g}H)$ and $c_{\sigma_g}(rH)=r\lr{g}H$ for all $r\in G$. 
Hence $\trans^G_\SSS$ satisfies the assumption of 
Theorem \ref{thm:abstract} with $\boldGamma=\trans^G_\SSS$. The abstract 
Burnside ring $\bz\trans^G_\SSS$ is isomorphic to the Burnside ring 
$\mathrm{\Omega}(G)$. 
\end{example}
\par
Let $R(G)$ be the ring of virtual $\bc$-characters, and let $R_\bq(G)$ be 
the subring of $R(G)$ generated by the characters afforded by $\bq G$-modules. 
There exists a ring homomorphism $\chara_G:\mathrm{\Omega}(G)\to R_\bq(G)$ 
given by 
\[
[G/H]\mapsto{1_H}^G 
\]
for all $H\leq G$, where ${1_H}^G$ is the character induced from the trivial 
character $1_H$ of $H$ (or the permutation character of $G$ on $G/H$) defined 
by 
\[
{1_H}^G(g)=\{rH\in G/H\st g\in\r{r}\!H\}
\]
for all $g\in G$. Obviously, if $u$ is a unit of $\mathrm{\Omega}(G)$, then 
$\chara_G(u)$ is a unit of $R_\bq(G)$. 
\par
For any unital ring $R$, we denote by $R^\times$ the unit group of $R$. 
Let $\hom(G,\lr{-1})$ be the set of homomorphisms 
from $G$ to the subgroup $\lr{-1}$ of $\bc^\times$. 
\begin{lemma}
\label{lem:idem}
For any $\chi\in R_\bq(G)^\times$, $\chi(\epsilon)\chi\in\hom(G,\lr{-1})$. 
\end{lemma}
{\itshape Proof.}
Obviously, $\chi(g)\in\lr{-1}$ for all $g\in G$, and 
the assertion follows from the first orthogonality relation 
(cf. \cite[(9.21), (9.26) Proposition]{CR87}). 
$\Box$
\par\bigskip
Let $I$ be an object of $\boldGamma$. We define an additive map 
$\omega_I:\bz\boldGamma\to\mathrm{\Omega}(\aut(I))$ by 
\[
[J]\mapsto[\boldGamma(I,J)] 
\]
for all objects $J$ of $\boldGamma$, where the left action of $\aut(I)$ on 
$\boldGamma(I,J)$ is given by 
\[
\sigma f=f\circ\sigma^{-1}
\]
for all $\sigma\in\aut(I)$ and $f\in\boldGamma(I,J)$. 
\par
There is a criteria for the units of an abstract Burnside ring, which is 
a generalization of that for the units of $\mathrm{\Omega}(G)$ 
(cf. \cite[Proposition 6.5]{Yo90a}). 
\begin{proposition}
\label{pro:Yoshida}
Keep the assumptions of Theorem {\rm\ref{thm:abstract}}, and assume further 
that for any object $I$ of $\boldGamma$, the map 
$\omega_I:\bz\boldGamma\to\mathrm{\Omega}(\aut(I))$ is a ring homomorphism. 
Let 
$\widetilde{x}=(x_I)_{[I]\in\boldGamma\!/\!\simeq}\in\bz^{\boldGamma\times}$. 
For any $[I]\in\boldGamma\!/\!\simeq$, define a map 
$\gamma_I^{\widetilde{x}}:\aut(I)\to\lr{-1}$ by 
\[
\sigma\mapsto x_Ix_{I/\sigma}
\]
for all $\sigma\in\aut(I)$. Then there exists an element $x$ of 
$\bz\boldGamma$ such that $\widetilde{x}=\varphi_\boldGamma(x)$, where 
$\varphi_\boldGamma=\varphi_\boldGamma^{(\infty)}$, if and only if 
$\gamma_I^{\widetilde{x}}\in\hom(\aut(I),\lr{-1})$ for any 
$[I]\in\boldGamma\!/\!\simeq$. 
\end{proposition}
{\itshape Proof.}
Let 
$\widetilde{x}=(x_I)_{[I]\in\boldGamma\!/\!\simeq}\in\bz^{\boldGamma\times}$. 
If $\gamma_I^{\widetilde{x}}\in\hom(\aut(I),\lr{-1})$ for any 
$[I]\in\boldGamma\!/\!\simeq$, then it follows from Theorem \ref{thm:abstract} 
that $\widetilde{x}=\varphi_\boldGamma(x)$ for some $x\in\bz\boldGamma$, 
because 
\[
\lr{\gamma_I^{\widetilde{x}},1_{\aut(I)}}_{\aut(I)}=
\dfrac{1}{\,|\aut(I)|\,}\sum_{\sigma\in\aut(I)}x_Ix_{I/\sigma}
\in\{0,\,1\}, 
\]
where $1_{\aut(I)}$ is the trivial character of $\aut(I)$ and 
$\lr{\gamma_I^{\widetilde{x}},1_{\aut(I)}}_{\aut(I)}$ is the inner product of 
$\gamma_I^{\widetilde{x}}$ and $1_{\aut(I)}$, for any 
$[I]\in\boldGamma\!/\!\simeq$. Conversely, assume that 
$\widetilde{x}=\varphi_\boldGamma(x)$ for some $x\in\bz\boldGamma$. Let $I$ 
be an object of $\boldGamma$. Since the map 
$\omega_I:\bz\boldGamma\to\mathrm{\Omega}(\aut(I))$ is a ring homomorphism, it 
follows that $\chara_{\aut(I)}(\omega_I(x))\in R_\bq(\aut(I))^\times$. By 
the definition of coequalizer, we have 
$x_{I/\sigma}=\chara_{\aut(I)}(\omega_I(x))(\sigma)$ for all 
$\sigma\in\aut(I)$ (cf. \cite[(2.28)]{YOT18}). Consequently, it follows from 
Lemma \ref{lem:idem} that $\gamma_I^{\widetilde{x}}\in\hom(\aut(I),\lr{-1})$. 
The proof is now complete. 
$\Box$
\par\bigskip
As before, we suppose that $\L$ is a finite $G$-lattice and 
$M_\L=(M_\L,\con{}{},\res{}{})$ is the monoid functor given in 
Proposition \ref{pro:lattice}. Let $\trans^G_\SSS$ be the category given in 
Example \ref{epl:trans}. There is an additive contravariant functor 
$\dot{T}^{M_\L}_G:\trans^G_\SSS\to\Mon$ inherited from 
$T^{M_\L}_G:G\set\to\Mon$. For each $K\leq G$, $T^{M_\L}_G(G/K)$ consists of 
the $G$-equivariant maps $\pi_s:G/K\to\widetilde{M_\L}(G),\,rK\mapsto\r{r}\!s$ 
for $s\in\L_K$. We call a pair $(G/K,\pi_s)$ of $G/K$ and $\pi_s$ with 
$(K,s)\in\SS(G,M_\L)$ an element of $\dot{T}^{M_\L}_G$. The morphisms 
$\lambda:(G/U,\pi_t)\to(G/K,\pi_s)$ between elements $(G/U,\pi_t)$ and 
$(G/K,\pi_s)$ of $\dot{T}^{M_\L}_G$ are defined to be the $G$-equivariant maps 
$\lambda_g:G/U\to G/K,\,rU\mapsto rgK$ for $gK\in\inv_U(G/K)$ 
(see Example \ref{epl:trans}) such that $t\leq\r{g}\!s$. Thus we obtain 
the category $\El(\dot{T}^{M_\L}_G)$ of elements of $\dot{T}^{M_\L}_G$. Note 
that the above definition of morphisms is different from that in the category 
of elements of $T^{M_\L}_G$. 
\par
For each $(U,t)\in\SS(G,M_\L)$, the group $\aut((G/U,\pi_t))$ of automorphisms 
of $(G/U,\pi_t)$ consists of all maps $\sigma_g:G/U\to G/U$ for 
$gU\in W_G(U,t)$ given by 
\[
rU\mapsto rgU
\]
for all $rU\in G/U$, which is identified with $\aut((G/U)_t)$. 
\par
We are now in a position to prove that $\bz_{(p)}\boldGamma$ with 
$\boldGamma=\El(\dot{T}^{M_\L}_G)$ is an abstract Burnside ring. By 
the following theorem, the lattice Burnside ring 
$\mathrm{\Omega}(G,M_\L)_{(p)}$, which is called a $p$-local lattice Burnside 
ring, is isomorphic to $\bz_{(p)}\boldGamma$. 
\begin{theorem}
\label{thm:trans}
Suppose that $\boldGamma=\El(\dot{T}^{M_\L}_G)$. Then $\boldGamma$ satisfies 
the assumptions of Theorem {\rm\ref{thm:abstract}}, and the abstract Burnside 
ring $\bz_{(p)}\boldGamma$ is isomorphic to the $p$-local lattice Burnside 
ring $\mathrm{\Omega}(G,M_\L)_{(p)}$. In this connection, 

\begin{equation}
\label{eq:coequalizer}
\begin{array}{l}
(G/U,\pi_t)/\sigma_g=(G/(\lr{g}U),\pi_{s_{(g,t)}}), 
\end{array}
\end{equation}
where $s_{(g,t)}=\inf\L_{\lr{g}U}^{\geq t}$, for all $(U,t)\in\SS(G,M_\L)$ and 
$gU\in W_G(U,t)$. Moreover, 
\begin{equation}
\label{eq:hom}
|\hom((G/U,\pi_t),(G/K,\pi_s))|=|\{gK\in G/K\st U\leq\r{g}\!K
\quad\mbox{and}\quad t\leq\r{g}\!s\}| 
\end{equation}
for all $(K,s),\,(U,t)\in\SS(G,M_\L)$ {\rm(}cf. Eq.{\rm(\ref{eq:varphi})}{\rm)}. 
\end{theorem}
{\itshape Proof.}
Since all the morphisms of $\boldGamma$ are epimorphisms, it follows that 
$\boldGamma$ has an epi-mono factorization system. Let $(U,t)\in\SS(G,M_\L)$, 
and let $g\in N_G(U,t)$. Observe that there exists a coequalizer 
$c_{\sigma_g}:(G/U,\pi_t)\to(G/U,\pi_t)/\sigma_g$ of $\sigma_g$ and 
$\id_{(G/U,\pi_t)}$ given by 
$(G/U,\pi_t)/\sigma_g=(G/(\lr{g}U),\pi_{s_{(g,t)}})$ and 
$c_{\sigma_g}(rU)=r\lr{g}U$ for all $r\in G$. Thus Eq.(\ref{eq:coequalizer}) 
holds, and $\boldGamma$ is a finite category which satisfies the assumptions 
of Theorem {\rm\ref{thm:abstract}}. Given $(K,s),\,(U,t)\in\SS(G,M_\L)$, 
$(G/K,\pi_s)\simeq(G/U,\pi_t)$ if and only if $(K,s)=g.(U,t)$ for some 
$g\in G$. Moreover, Eq.(\ref{eq:hom}) is obvious. Consequently, it follows 
from Proposition \ref{pro:monoid} and Theorems \ref{thm:fundamental} and 
\ref{thm:ring homomorphism} that there is a ring isomorphism 
$\mathrm{\Omega}(G,M_\L)_{(p)}\stackrel{\sim}{\to}\bz_{(p)}\boldGamma$ given 
by 
\[
[(G/K)_s]\mapsto[(G/K,\pi_s)]
\]
for all $(K,s)\in\SS(G,M_\L)$. This completes the proof. 
$\Box$
\par\bigskip
While the ring structure of $\mathrm{\Omega}(G,M_\L)$ is given in 
Proposition \ref{pro:monoid}, the proof of the following corollary to 
Theorem \ref{thm:trans} makes it clear by a Burnside map. 
\begin{corollary}
\label{cor:ring homomorphism}
Suppose that $\boldGamma=\El(\dot{T}^{M_\L}_G)$, and let $I$ be an object of 
$\boldGamma$. Then the map $\omega_I:\bz\boldGamma\to\mathrm{\Omega}(\aut(I))$ 
is a ring homomorphism. 
\end{corollary}
{\itshape Proof.}
Let $(U,t), (K_i,s_i)\in\SS(G,M_\L)$ with $i=1,\,2$, and define a map 
\[
\begin{array}{l}
\vspace{0.3cm}
\Phi:\boldGamma((G/U,\pi_t),(G/K_1,\pi_{s_1}))\times
\boldGamma((G/U,\pi_t),(G/K_2,\pi_{s_2}))\\
\hspace{3cm}\displaystyle\to\dot{\displaystyle
\bigcup_{K_1rK_2\in K_1\backslash G/K_2}}
\boldGamma((G/U,\pi_t),(G/(K_1\cap\r{r}\!K_2),\pi_{s_1\wedge\r{r}\!s_2}))
\end{array}
\]
by 
\[
(\lambda_{g_1},\lambda_{g_2})\mapsto\lambda_{g_1r_1}
\in\boldGamma((G/U,\pi_t),(G/(K_1\cap\r{r}\!K_2),\pi_{s_1\wedge\r{r}\!s_2})), 
\]
where $g_1^{-1}g_2K_2=r_1rK_2$ with $r_1\in K_1$, for all 
$g_iK_i\in\inv_U(G/K_i)$ with $i=1,\,2$ such that $t\leq\r{g_i}\!s_i$. Then $\Phi$ is an isomorphism of 
$\aut((G/U,\pi_t))$-sets. Hence it follows from Theorem \ref{thm:trans} that 
$\omega_I$ is a ring homomorphism. 
$\Box$
\par\bigskip
In the proof of Theorem \ref{thm:trans}, we are aware that the exact sequence 
{\rm(\ref{eq:abstract})} with $\boldGamma=\El(\dot{T}^{M_\L}_G)$ is identified 
with the exact sequence {\rm(\ref{eq:lattice})} with $H=G$, namely, 
\[
0\longrightarrow\mathrm{\Omega}(G,M_\L)_{(p)}
\overset{\alpha_G^{(p)}\circ\varphi_G^{(p)}}
{\longrightarrow}\widetilde{\mathrm{\Omega}}(G,M_\L)_{(p)}
\overset{\widetilde{\psi}_G^{(p)}\circ\beta_G^{(p)}}{\longrightarrow}
\obs(G,M_\L)_{(p)}\longrightarrow0. 
\]
(By Theorems \ref{thm:fundamental} and \ref{thm:ring homomorphism} with $H=G$, 
$\widetilde{\alpha}_G^{(p)}=\alpha_G^{(p)}\circ\varphi_G^{(p)}$, 
$\widetilde{\beta}_G^{(p)}=\widetilde{\psi}_G^{(p)}\circ\beta_G^{(p)}$, and 
the map $\alpha_G\circ\varphi_G$ is a ring monomorphism, where 
$\alpha_G=\alpha_G^{(\infty)}$.) 
\par
We are now ready to give a criteria of the units of $\mathrm{\Omega}(G,M_\L)$. 
\begin{proposition}
\label{pro:Yoshida criteria}
Let $\widetilde{x}=(x_{(U,t)})_{(U,t)\in\R(G,M_\L)}\in
\widetilde{\mathrm{\Omega}}(G,M_\L)^\times$. Then $\widetilde{x}$ is contained 
in the image of the ring monomorphism 
$\alpha_G\circ\varphi_G:\mathrm{\Omega}(G,M_\L)\to
\widetilde{\mathrm{\Omega}}(G,M_\L)$ if and only if, for each 
$(U,t)\in\R(G,M_\L)$, the map 
$\gamma_{(U,t)}^{\widetilde{x}}:\aut((G/U)_t)\to\lr{-1}$ given by 
\[
\sigma_g\mapsto x_{(U,t)}x_{(\lr{g}U,s_{(g,t)})} 
\]
for all $gU\in W_G(U,t)$ is a group homomorphism. 
\end{proposition}
{\itshape Proof.}
The proposition follows from Proposition \ref{pro:Yoshida} and 
Corollary \ref{cor:ring homomorphism}. 
$\Box$
\par\bigskip
Let $M_\SSS$ be the monoid functor given in Example \ref{epl:slice}. Then 
Proposition \ref{pro:Yoshida criteria} with $M_\L=M_\SSS$ is equivalent to 
\cite[Theorem 8.4]{Bou12}. 
\par
We give an example of $\mathrm{\Omega}(G,M_\L)^\times$ (see also 
\cite[Theorem A.13]{Bou12}). 
\begin{proposition}
\label{pro:unit lattice}
Suppose that $G$ is abelian. Let $M_\SSS^\circ$ be the monoid functor given in 
Example {\rm\ref{epl:coslice}}. Then $\mathrm{\Omega}(G,M_\SSS^\circ)^\times$ 
is generated by $-[(G/G)_G]$ and the elements $[(G/H)_U]-[(G/G)_G]$ for 
$U\leq H\leq G$ with $|G:H|=2$. In particular, 
$\mathrm{\Omega}(G,M_\SSS^\circ)^\times$ is an elementary abelian $2$-group of 
order $2^{\vartheta(G)+1}$, where 
\[
\vartheta(G)=|\{(H,U)\in\SS(G,M_\SSS^\circ)\st|G:H|=2\}|. 
\]
\end{proposition}
{\itshape Proof.}
Let $\widetilde{x}=(x_{(K,L)})_{(K,L)\in\R(G,M_\SSS^\circ)}\in
\widetilde{\mathrm{\Omega}}(G,M_\SSS^\circ)^\times$, and suppose that 
$\widetilde{x}$ is contained in the image of the ring monomorphism 
$\alpha_G\circ\varphi_G:\mathrm{\Omega}(G,M_\L)\to
\widetilde{\mathrm{\Omega}}(G,M_\L)$ with $M_\L=M_\SSS^\circ$. Then by 
Proposition \ref{pro:Yoshida criteria}, we have 
\[
x_{(K,L)}x_{(\lr{g}K,L)}x_{(\lr{r}K,L)}=x_{(\lr{gr}K,L)}
\]
for all $(K,L)\in\R(G,M_\SSS^\circ)(=\SS(G,M_\SSS^\circ))$ and $gK,\,rK\in W_G(K,L)$. This 
implies that for each $(K,L)\in\R(G,M_\SSS^\circ)$, if $|G:K|>2$, then 
the value $x_{(K,L)}$ is determined by the values $x_{(H,L)}$ for 
$H\leq G$ with $K<H$ (cf. \cite[p. 904]{Bou12}). Hence we have 
\begin{equation}
\label{eq:order}
|\mathrm{\Omega}(G,M_\SSS^\circ)^\times|\leq2^{\vartheta(G)+1}. 
\end{equation}
If $G$ is of odd order, then the assertion clearly holds. Assume that $G$ is 
of even order, and let $H_1,\,H_2,\dots,\,H_m$ be the subgroups of index $2$ 
in $G$. For each integer $i$ with $1\leq i\leq m$, let $\U(H_i)$ denote 
the subgroup of $\mathrm{\Omega}(G,M_\SSS^\circ)^\times$ generated by 
the elements $[(G/H_i)_U]-[(G/G)_G]$ for $U\leq H_i$. By 
Eq.(\ref{eq:product}), we have 
\[
(\U(H_1)\cdots\U(H_i))\cap\U(H_{i+1})=\{[(G/G)_G]\} 
\]
for each integer $i$ with $1\leq i\leq m-1$ and 
\[
(\U(H_1)\cdots\U(H_m))\cap\lr{-[(G/G)_G]}=\{[(G/G)_G]\}. 
\]
Thus $\lr{-[(G/G)_G]}\U(H_1)\cdots\U(H_m)=\lr{-[(G/G)_G]}\times\U(H_1)\times
\cdots\times\U(H_m)$. Likewise, for each integer $i$ with $1\leq i\leq m$, 
$\U(H_i)$ is the direct product of the subgroups $\lr{[(G/H_i)_U]-[(G/G)_G]}$ 
for $U\leq H_i$. By these facts, $\mathrm{\Omega}(G,M_\SSS^\circ)^\times$ 
contains the direct product of the subgroups $\lr{[(G/H_i)_U]-[(G/G)_G]}$ for 
$i=1,\,2,\dots,\,m$ and $U\leq H_i$. Combining this fact with 
Eq.(\ref{eq:order}), we conclude that the assertion holds. 
$\Box$
\begin{remark}
\label{remark}
There is an embedding $\mathrm{\Omega}(G)\hookrightarrow
\mathrm{\Omega}(G,M_\SSS^\circ)$ given by 
\[
[G/H]\mapsto[(G/H)_H]
\]
for all $H\leq G$. By Proposition \ref{pro:unit lattice}, 
$\mathrm{\Omega}(G)^\times $ is generated by $-[G/G]$ and the elements 
$[G/H]-[G/G]$ for $H\leq G$ with $|G:H|=2$, and 
$|\mathrm{\Omega}(G)^\times|=2^{|\Hom(G,\lr{-1})|}$ (see also 
\cite[Remark A.14]{Bou12} and \cite[Lemma 7.1]{Yo90a}), which is due to 
Matsuda \cite[Example 4.5]{Mat82}. 
\end{remark}
%----------------------------------------------------------------
%
%newsection
%
%----------------------------------------------------------------
\section{Primitive idempotents of lattice Burnside rings}
%----------------------------------------------------------------
Let $\boldGamma$ be a finite category, and suppose that the assumptions of 
Theorem \ref{thm:abstract} hold. Let 
$\id_\bq\otimes\varphi_\boldGamma:\bq\otimes_\bz\bz\boldGamma\to
\bq\otimes_\bz\bz^\boldGamma$ be the algebra monomorphism determined by 
$\varphi_\boldGamma^{(\infty)}$. By Theorem \ref{thm:abstract}, the primitive 
idempotents of $\bq\otimes_\bz\bz\boldGamma$ are the elements $e_I$ for 
$[I]\in\boldGamma\!/\!\simeq$ such that $\id_\bq\otimes\varphi_\boldGamma(e_I)
=(\delta_{[I]\,[J]})_{[J]\in\boldGamma\!/\!\simeq}$. 
\par
Let $\sim_p$ be the equivalence relation on the set $\boldGamma\!/\!\simeq$ 
generated by 
\[
[I/\sigma]\sim_p[I]\quad\mbox{with}\quad\sigma\in\aut(I)_p. 
\]
(Note that $[I]=[I/\id_I]\sim_p[I]$ for any object $I$ of $\boldGamma$.) 
We define an equivalence relation $\sim_p$ on the set of objects of 
$\boldGamma$ by letting 
\[
I\sim_pJ\quad\mbox{if and only if}\quad[I]\sim_p[J]. 
\]
\par
Let $\C_p(\boldGamma)$ be a complete set of representatives of equivalence 
classes with respect to the equivalence relation $\sim_p$ on $\boldGamma$. For 
each $I\in\C_p(\boldGamma)$, we define
\[
e_I^{(p)}:=\sum_{[I]\sim_p[J]\in\boldGamma\!/\!\simeq}e_J, 
\]
where the sum is taken over all $[J]\in\boldGamma\!/\!\simeq$ such that 
$I\sim_pJ$. 
\par
There is a generalization of \cite[Lemma 2]{Glu81} and 
\cite[Theorem 3.1]{Yo83a}: 
\begin{proposition}
\label{pro:idem}
The primitive idempotents of $\bz_{(p)}\boldGamma$ coincide with the elements 
$e_I^{(p)}$ for $I\in\C_p(\boldGamma)$, and those of $\bz\boldGamma$ coincide 
with the elements $e_I^{(\infty)}$ for $I\in\C_\infty(\boldGamma)$. 
\end{proposition}
{\itshape Proof.}
Let $(x_I)_{[I]\in\boldGamma\!/\!\simeq}$ be an idempotent of 
$\bz_{(p)}^\boldGamma$. Then by Theorem \ref{thm:abstract}, 
$(x_I)_{[I]\in\boldGamma\!/\!\simeq}$ is contained in the image of 
$\varphi_\boldGamma^{(p)}$ if and only if $x_I=x_J=0$ or $x_I=x_J=1$ for all 
pairs $(I,J)$ with $I\sim_pJ$. Consequently, the primitive idempotents of  
$\bz_{(p)}\boldGamma$ coincide with the elements $e_I^{(p)}$ for 
$I\in\C_p(\boldGamma)$, as desired.
$\Box$
\par\bigskip
We turn to the primitive idempotents of $\mathrm{\Omega}(G,M_\L)$, where $\L$ 
is a finite $G$-lattice and $M_\L=(M_\L,\con{}{},\res{}{})$ is the monoid 
functor given in Proposition \ref{pro:lattice}. 
\par
Set $\boldGamma=\dot{T}^{M_\L}_G$. By Theorem \ref{thm:trans}, $\boldGamma$ 
satisfies the assumptions of Theorem {\rm\ref{thm:abstract}}, and the abstract 
Burnside ring $\bz_{(p)}\boldGamma$ is isomorphic to the $p$-local lattice 
Burnside ring $\mathrm{\Omega}(G,M_\L)_{(p)}$. We define an equivalence 
relation $\sim_p$ on $\SS(G,M_\L)$ by letting 
\[
(K,s)\sim_p(U,t)\quad\mbox{if and only if}\quad(G/K,\pi_s)\sim_p(G/U,\pi_t). 
\]
Given $(U,t)\in\SS(G,M_\L)$ and $gU\in W_G(U,t)$, it follows from 
Theorem \ref{thm:trans} that 
\[
(G/(\lr{g}U),\pi_{s_{(g,t)}})=(G/U,\pi_t)/\sigma_g\sim_p(G/U,\pi_t)
\]
with $s_{(g,t)}=\inf\L_{\lr{g}U}^{\geq t}$, whence 
$(\lr{g}U,s_{(g,t)})\sim_p(U,t)$. We often use this basic fact. 
\par
Let $K\leq G$. When $p$ is a prime, we denote by $O^p(K)$ the smallest normal 
subgroup of $K$ such that $K/O^p(K)$ is a $p$-group. Suppose that 
\[
K=K^{(0)}\geq K^{(1)}\geq K^{(2)}\geq\cdots\geq K^{(i)}\geq\cdots
\]
is the derived series of $K$ (cf. \cite[Chapter 2, Definition 3.11]{Su82}). 
Then we define $O^\infty(K):=\cap_{i=1}^\infty K^{(i)}$. A subgroup $K$ of $G$ 
is said to be $p$-perfect if $K=O^p(K)$. 
\begin{lemma}
\label{lem:der}
Let $(K,s),\,(U,t)\in\SS(G,M_\L)$, and assume that $(K,s)\sim_p(U,t)$. Then 
the subgroup $O^p(K)$ of $K$ is a conjugate of the subgroup $O^p(U)$ of $U$ in 
$G$. 
\end{lemma}
{\itshape Proof.}
We may assume that $K=\lr{g}U$ for some $gU\in N_G(U,t)_p$. If $p$ is a prime, 
then $U\geq O^p(U)\geq O^p(K)$, and thus $O^p(U)=O^p(K)$. Suppose that 
$p=\infty$. We have $U^{(i-1)}\geq K^{(i)}\geq U^{(i)}$ for any $i\geq1$. If 
$U^{(i-1)}=U^{(i)}$ for some $i$, then $U^{(i-1)}=K^{(i)}=U^{(i)}$. Hence we 
have $O^\infty(K)=O^\infty(U)$, completing the proof. 
$\Box$
\begin{lemma}
\label{lem:inf}
Let $U$ and $H$ be subgroups of $G$, and suppose that $O^p(U)$ is a conjugate 
of $O^p(H)$ in $G$. Then $(H,\inf\L_H)\sim_p(U,\inf\L_U)$. 
\end{lemma}
{\itshape Proof.}
We may assume that $H=\lr{g}U$ for some $gU\in W_G(U,\inf\L_U)_p$. By 
definition, $\inf\L_{\lr{g}U}\wedge\sup\L_U\in\L_U$. Set $t=\inf\L_U$. Then 
$t\leq\inf\L_{\lr{g}U}=\inf\L_{\lr{g}U}^{\geq t}$. Hence we 
have $(\lr{g}U,\inf\L_{\lr{g}U})\sim_p(U,\inf\L_U)$, completing the proof. 
$\Box$
\par\bigskip
Given $K\leq G$ and $s_1,\,s_2\in\L_K$ with $s_1>s_2$, the phrase 
`$s_1$ covers $s_2$ in $\L_K$' means that there is no element $t$ of $\L_K$ 
satisfying the condition $s_1>t>s_2$. 
\par
We now extend Dress' characterization of solvable groups for 
$\mathrm{\Omega}(G)$ (cf. \cite{Dr69}). 
\begin{theorem}
\label{thm:solvable}
Assume that, given $s\in\L$ with $s\in\L_K$ for some $K\leq G$, there exist 
subgroups $H$ and $U$ of $G$ with $U\leq H\leq N_G(U)$ satisfying 
the following conditions\! {\rm:} 
\def\theenumi{\roman{enumi}}
\begin{enumerate}
\item
\label{thm:solvable1}
The set of subgroups $K$ of $H$ containing $U$ coincides with 
$\{K\leq G\st s\in\L_K\}$. 
\item
\label{thm:solvable2}
If $s$ covers $t$ in $\L_U$, then $s=\inf\L_{\lr{g}U}^{\geq t}$ for some 
$gU\in W_G(U,t)_p$ with $g\not\in U$. 
\end{enumerate}
\def\theenumi{\alph{enumi}}
Then $G$ is a $p$-group, where an $\infty$-group is a solvable group, if and 
only if the prime spectrum of $\mathrm{\Omega}(G,M_\L)_{(p)}$ is connected in 
the Zariski topology, that is, if and only if $0$ and $1$ are the only 
idempotents of $\mathrm{\Omega}(G,M_\L)_{(p)}$. 
\end{theorem}
{\itshape Proof.}
By Proposition \ref{pro:idem}, $0$ and $1$ are the only idempotents of 
$\mathrm{\Omega}(G,M_\L)_{(p)}$ if and only if $(K,s)\sim_p(U,t)$ for all 
$(K,s),\,(U,t)\in\SS(G,M_\L)$. If $\SS(G,M_\L)$ has only one equivalence class 
with respect to $\sim_p$, then by Lemma \ref{lem:der}, $O^p(G)=\{\epsilon\}$, 
which forces $G$ to be a $p$-group. Conversely, assume that $G$ is 
a $p$-group. Let $K\leq G$, and let $s\in\L_K$. By 
the condition (\ref{thm:solvable1}), there exist subgroups $H$ and $U$ of $G$ 
such that $U\leq H\leq N_G(U,s)$ and $\{K_1\leq G\st s\in\L_{K_1}\}$ coincides 
with the set of subgroups $K_1$ of $H$ containing $U$. We have 
$(K,s)\sim_p(U,s)$, because $K/U$ is a $p$-group. If $s$ covers $t$ in $\L_U$, 
then by the condition (\ref{thm:solvable2}), $s=\inf\L_{\lr{g}U}^{\geq t}$ and 
$(\lr{g}U,s)\sim_p(U,t)$ for some $gU\in W_G(U,t)_p$ with $g\not\in U$, which 
implies that $(U,s)\sim_p(\lr{g}U,s)\sim_p(U,t)$. Hence either $s=\inf\L_U$ or 
$(K,s)\sim_p(U,t)$ for some $t\in\L_U$ with $t<s$. By repeating this argument, 
we can choose elements $(U_0,t_0):=(K,s),\,(U_1,t_1),\dots,\,(U_\ell,t_\ell)$ 
of $\SS(G,M_\L)$ with $t_\ell=\inf\L_{U_\ell}$ such that 
$(U_i,t_i)\sim_p(U_{i+1},t_{i+1})$ and $t_{i+1}\leq t_i\in\L_{U_{i+1}}$ with 
$i=0,\,1,\dots,\,\ell-1$. Moreover, 
$(U_\ell,\inf\L_{U_\ell})\sim_p(\{\epsilon\},\inf\L_{\{\epsilon\}})$ by 
Lemma \ref{lem:inf}. Thus we have $(K,s)\sim_p
(U_\ell,\inf\L_{U_\ell})\sim_p(\{\epsilon\},\inf\L_{\{\epsilon\}})$. 
Consequently, $\SS(G,M_\L)$ has only one equivalence class with respect to 
$\sim_p$. This completes the proof. 
$\Box$
\begin{example}
\label{epl:slice solvable}
Keep the notation of Example \ref{epl:slice}, and assume further that 
$G$ is a $p$-group. Let $E\leq G$. Then $\{K\leq G\st E\in\SSS_{\geq K}\}$ is 
the set of subgroups of $E$. Suppose that $E\not=\{\epsilon\}$ and $E$ covers 
$F$ in $\SSS_{\geq\{\epsilon\}}$. Since $G$ is a $p$-group, it turns out that 
$F$ is a normal maximal subgroup of $E$. We can take an element $g$ of 
$N_E(F)$ for which $E=\lr{g}F$. Then it is obvious that 
$E=\inf\SSS_{\geq\lr{g}F}$. Hence the assumption of Theorem \ref{thm:solvable} 
with $M_\L=M_\SSS$ and $s=E\in\SSS$ holds. Consequently, $0$ and $1$ are 
the only idempotents of $\mathrm{\Omega}(G,M_\SSS)_{(p)}$. The assertion of 
Theorem \ref{thm:solvable} with $M_\L=M_\SSS$ is given in 
\cite[Theorem 7.9]{Bou12}, together with the primitive idempotents of 
$\mathrm{\Omega}(G,M_\SSS)$. 
\end{example}
\par
This section ends with a deference between $\mathrm{\Omega}(G,M_\SSS)_{(p)}$ 
and $\mathrm{\Omega}(G,M_\SSS^\circ)_{(p)}$. 
\begin{theorem}
\label{thm:normal}
Let $M_\SSS^\circ$ be the monoid functor given in 
Example {\rm\ref{epl:coslice}}. Then $G$ is a $p$-group if and only if 
the number of primitive idempotents of $\mathrm{\Omega}(G,M_\SSS^\circ)_{(p)}$ 
is $|\C(G)|$, where $\C(G)$ is a full set of nonconjugate subgroups of $G$. 
\end{theorem}
{\itshape Proof.}
Let $(H,U)\in\SS(G,M_\SSS^\circ)$, and let $gU\in W_G(H,U)_p$. Since 
$\lr{g}H\leq N_G(U)$, we have $(\lr{g}H,U)\in\SS(G,M_\SSS^\circ)$ and 
$U=\inf\SSS_{\unlhd\lr{g}H}^{\geq U}$. Hence, if $(K,L)\sim_p(H,U)$ for some 
$(K,L)\in\SS(G,M_\SSS^\circ)$, then $L$ is a conjugate of $U$ in $G$. Assume 
that the number of primitive idempotents of 
$\mathrm{\Omega}(G,M_\SSS^\circ)_{(p)}$ is $|\C(G)|$. By 
Proposition \ref{pro:idem} and the above fact, we have 
$(G,\{\epsilon\})\sim_p(\{\epsilon\},\{\epsilon\})$. Thus it follows from 
Lemma \ref{lem:der} that $O^p(G)=\{\epsilon\}$, which forces $G$ to be 
a $p$-group. Conversely, assume that $G$ is a $p$-group. Then by the previous 
fact, $(H,U)\sim_p(U,U)$ for any $U\leq H\leq N_G(U)$. Moreover, for any 
subgroups $K$ and $U$ of $G$, if $(K,K)\sim_p(U,U)$, then $K$ is a conjugate 
of $U$ in $G$. Consequently, it follows from Proposition \ref{pro:idem} that 
the number of primitive idempotents of $\mathrm{\Omega}(G,M_\SSS^\circ)_{(p)}$ 
is $|\C(G)|$. This completes the proof. 
$\Box$
%----------------------------------------------------------------
%
%newsection
%
%----------------------------------------------------------------
\section{Partial lattice Burnside rings }
%----------------------------------------------------------------
We define a partially order $\leq$ on $\SS(G,M_\L)$ by the rule that 
\[
(U,t)\leq(K,s)\quad\Longleftrightarrow\quad U\leq K\quad\mbox{and}\quad 
t\leq s. 
\]
\par
Let $\XX$ be a subset of $\SS(G,M_\L)$ closed under the action of $G$ (see 
Definition \ref{def:pi}), and suppose that the following condition holds. 
\begin{condition}
Given $(U,t)\in\XX$ and $gU\in W_G(U,t)$, the set 
\[
\{(K,s)\in\XX\st(K,s)\geq(\lr{g}U,s_{(g,t)})\}, 
\]
where $s_{(g,t)}=\inf\L_{\lr{g}U}^{\geq t}$, has a unique minimal element, 
denoted by $(\ov{\lr{g}U},\ov{s_{(g,t)}})$, with respect to the partially 
order $\leq$ on $\SS(G,M_\L)$. 
\end{condition}
\def\theenumi{\alph{enumi}}
\par
We denote by $\El(\dot{T}^{M_\L}_G)_\XX$ the full subcategory of 
$\El(\dot{T}^{M_\L}_G)$ whose objects are the elements $(G/K,\pi_s)$ for 
$(K,s)\in\XX$, and write $\boldGamma_\XX=\El(\dot{T}^{M_\L}_G)_\XX$. 
\begin{lemma}
\label{lem:partial}
The category $\boldGamma_\XX$ is a finite category which satisfies 
the assumptions of Theorem {\rm\ref{thm:abstract}}, so that 
$\bz_{(p)}\boldGamma_\XX$ is an abstract Burnside ring. Moreover, 
\[
(G/U,\pi_t)/\sigma_g=(G/(\ov{\lr{g}U}),\pi_{\ov{s_{(g,t)}}}) 
\]
for all $(U,t)\in\XX$ and $gU\in W_G(U,t)$. 
\end{lemma}
{\itshape Proof.}
Let $(U,t)\in\XX$, and let $gU\in W_G(U,t)$. By the condition (A), there 
exists a coequalizer $c_{\sigma_g}:(G/U,\pi_t)\to(G/U,\pi_t)/\sigma_g$ of 
$\sigma_g$ and $\id_{(G/U,\pi_t)}$ given by 
$(G/U,\pi_t)/\sigma_g=(G/(\ov{\lr{g}U}),\pi_{s_{\ov{(g,t)}}})$ and 
$c_{\sigma_g}(rU)=r\lr{g}U$ for all $r\in G$, as desired. 
$\Box$
\par\bigskip
Set $\ov{\XX}=\R(G,M_\L)\cap\XX$. We denote by $\mathrm{\Omega}(G,\XX)_{(p)}$ 
the $\bz_{(p)}$-submodule of $\mathrm{\Omega}(G,M_\L)_{(p)}$ generated by 
the elements $[(G/K)_s]$ for $(K,s)\in\ov{\XX}$, and define 
\[
\widetilde{\mathrm{\Omega}}(G,\XX)_{(p)}:=\prod_{(K,s)\in\ov{\XX}}\bz_{(p)}
\quad\mbox{and}\quad
\obs(G,\XX)_{(p)}:=\prod_{(K,s)\in\ov{\XX}}\bz_{(p)}/|W_G(K,s)_p|\bz_{(p)}. 
\]
The $\bz_{(p)}$-module $\mathrm{\Omega}(G,\XX)_{(p)}$  is identified with 
the abstract Burnside ring $\bz_{(p)}\boldGamma_\XX$. We define 
a $\bz_{(p)}$-module homomorphism 
$\varphi_\XX^{(p)}:\mathrm{\Omega}(G,\XX)_{(p)}\to
\widetilde{\mathrm{\Omega}}(G,\XX)_{(p)}$ by 
\[
[(G/H)_s]\mapsto(|\{gH\in G/H\st U\leq\r{g}\!H\quad\mbox{and}\quad
t\leq\r{g}\!s\}|)_{(U,t)\in\ov{\XX}} 
\]
for all $(H,s)\in\XX$, and define a map 
$\psi_\XX^{(p)}:\widetilde{\mathrm{\Omega}}(G,\XX)_{(p)}\to\obs(G,\XX)_{(p)}$ 
by 
\[
(x_{(H,s)})_{(H,s)\in\ov{\XX}}\mapsto\left(\sum_{gU\in W_G(U,t)_p}
x_{(\ov{\lr{g}U},\ov{s_{(g,t)}})}\bmod|W_G(U,t)_p|\right)_{(U,t)\in\ov{\XX}}
\]
for all 
$(x_{(H,s)})_{(H,s)\in\ov{\XX}}\in\widetilde{\mathrm{\Omega}}(G,\XX)_{(p)}$; 
these maps are identified with the Burnside map and the Cauchy-Frobenius map, 
respectively (see Theorem \ref{thm:trans} and Lemma \ref{lem:partial}). 
\begin{theorem}
\label{thm:partial p}
Under the above notation, the sequence 
\[
0\longrightarrow\mathrm{\Omega}(G,\XX)_{(p)}\overset{\varphi_\XX^{(p)}}
{\longrightarrow}\widetilde{\mathrm{\Omega}}(G,\XX)_{(p)}
\overset{\psi_\XX^{(p)}}{\longrightarrow}\obs(G,\XX)_{(p)}\longrightarrow0
\]
of $\bz_{(p)}$-modules is exact. Moreover, the $\bz_{(p)}$-module 
$\mathrm{\Omega}(G,\XX)_{(p)}$ has a $\bz_{(p)}$-algebra structure such that 
$\varphi_\XX^{(p)}$ is an injective $\bz_{(p)}$-algebra homomorphism. 
\end{theorem}
{\itshape Proof.}
The theorem follows from Theorem \ref{thm:abstract} and 
Lemma \ref{lem:partial}. 
$\Box$
\begin{corollary}
\label{cor:section}
Let $\mu_\XX$ denote the M\"{o}bius function of the partially ordered set 
$\XX$ with the binary relation $\leq$. The primitive idempotents of 
$\bq\otimes_\bz\mathrm{\Omega}(G,\XX)$ are the elements 
\[
\varepsilon_{(K,s)}:=\frac{1}{\,|N_G(K,s)|\,}
\sum_{\genfrac{}{}{0pt}{3}{(U,t)\in\XX,}{(U,t)\leq(K,s)}}|U|
\mu_\XX((U,t),(K,s))[(G/U)_t] 
\]
for $(K,s)\in\ov{\XX}$, where the sum is taken over all $(U,t)\in\XX$ with 
$(U,t)\leq(K,s)$. 
\end{corollary}
{\itshape Proof.}
Set $\varphi_\XX=\varphi_\XX^{(\infty)}$. For each $(K,s)\in\XX$, we have 
\[
\begin{array}{l}
\vspace{0.1cm}
\varphi_\XX(|G|\cdot\varepsilon_{(K,s)})\\
\vspace{0.1cm}\hspace{1cm}\displaystyle
=\frac{|G|}{\,|N_G(K,s)|\,}
\sum_{\genfrac{}{}{0pt}{3}{(U,t)\in\XX,}{(U,t)\leq(K,s)}}|U|
\mu_\XX((U,t),(K,s))\varphi_\XX([(G/U)_t])\\
\vspace{0.1cm}\hspace{1cm}\displaystyle
=\frac{|G|}{\,|N_G(K,s)|\,}
\sum_{\genfrac{}{}{0pt}{3}{(U,t)\in\XX,}{(U,t)\leq(K,s)}}|U|
\mu_\XX((U,t),(K,s))\\
\vspace{0.1cm}\hspace{5cm}\displaystyle
\cdot(|\{gU\in G/U\st U_1\leq\r{g}U\quad\mbox{and}\quad t_1\leq\r{g}t\}|)
_{(U_1,t_1)\in\ov{\XX}}\\
\vspace{0.1cm}\hspace{1cm}\displaystyle
=\frac{|G|}{\,|N_G(K,s)|\,}\sum_{g\in G}\left
(\sum_{(U_,t)\in\XX_{(\r{g^{-1}}\!U_1,\r{g^{-1}}\!t_1)}^{(K,s)}}
\mu_\XX((U,t),(K,s))\right)_{(U_1,t_1)\in\ov{\XX}}\\
\hspace{1cm}\displaystyle
=|G|\cdot(\delta_{(K,s),(U,t)})_{(U,t)\in\ov{\XX}}, 
\end{array}
\]
where $\XX_{(\r{g^{-1}}\!U_1,\r{g^{-1}}\!t_1)}^{(K,s)}
=\{(U,t)\in\XX\st(U_1,t_1)\leq(\r{g}U,\r{g}t)\leq(\r{g}\!K,\r{g}\!s)\}$. Hence 
the assertion is an immediate consequence of Theorem {\rm\ref{thm:partial p}}. 
$\Box$
\begin{remark}
\label{rem:Mebius}
From Theorem \ref{thm:ring homomorphism} with $H=G$ and 
Corollary \ref{cor:section} with $\XX=\SS(G,M_\L)$, we know that the elements 
$e_{(K,s)}^{(G)}$ of $\bq\otimes_\bz\mathrm{\Omega}(G,M_\L)$ for 
$(K,s)\in\R(G,M_\L)$ are the primitive idempotents, and so are the elements 
$\varepsilon_{(K,s)}$ of $\bq\otimes_\bz\mathrm{\Omega}(G,M_\L)$ for 
$(K,s)\in\R(G,M_\L)$. Let $(K,s)\in\R(G,M_\L)$. Given $(U,t)\in\SS(G,M_\L)$, 
we define 
\[
\widetilde{\mu}((U,t),(K,s)):=
\left\{
\begin{array}{cl}
\vspace{0.1cm}\displaystyle
\sum_{\genfrac{}{}{0pt}{3}{t_1\in\L_K^{\leq s},}{t=t_1\wedge\sup\L_U}}
\mu_K(t_1,s)\mu(U,K)\quad&\mbox{if }(U,t)\leq(K,s), \\
0\quad&\mbox{otherwise}, 
\end{array}\right.
\]
where the sum is taken over all $t_1\in\L_K^{\leq s}$ with 
$t=t_1\wedge\sup\L_U$. Observe that 
\[
\varepsilon_{(K,s)}=\frac{1}{\,|N_G(K,s)|\,}
\sum_{\genfrac{}{}{0pt}{3}{(U,t)\in\XX,}{(U,t)\leq(K,s)}}|U|
\mu_\XX((U,t),(K,s))[(G/U)_t] 
\]
and
\[
e_{(K,s)}^{(G)}=\frac{1}{\,|N_G(K,s)|\,}
\sum_{\genfrac{}{}{0pt}{3}{(U,t)\in\SS(G,M_\L),}{(U,t)\leq(K,s)}}|U|
\widetilde{\mu}((U,t),(K,s))[(G/U)_t]. 
\]
Set $\XX=\SS(G,M_\L)$. We have 
$\mu_\XX((U,t),(K,s))=\widetilde{\mu}((U,t),(K,s))$, or equivalently, 
\[
\sum_{(U,t)\leq(H,s_1)\leq(K,s)}\widetilde{\mu}((U,t),(H,s_1))
=\left\{
\begin{array}{cl}
\vspace{0.1cm}\displaystyle
1\quad&\mbox{if }(U,t)=(K,s),\\
0\quad&\mbox{otherwise}, 
\end{array}
\right. 
\]
for all $(U,t)\in\XX$, by which $\varepsilon_{(K,s)}$ coincides with 
$e_{(K,s)}^{(G)}$. In fact, it follows that 
\[
\begin{array}{ll}
\vspace{0.1cm}\displaystyle
\sum_{(U,t)\leq(H,s_1)\leq(K,s)}\widetilde{\mu}((U,t),(H,s_1))&\displaystyle
=\sum_{(U,t)\leq(H,s_1)\leq(K,s)}
\sum_{\genfrac{}{}{0pt}{3}{t_1\in\L_H^{\leq s_1},}{t=t_1\wedge\sup\L_U}}
\mu_H(t_1,s_1)\mu(U,H)\\
\vspace{0.1cm}
\mbox{}&\displaystyle
=\sum_{U\leq H\leq K}
\sum_{\genfrac{}{}{0pt}{3}{t_1\in\L_H^{\leq s},}{t=t_1\wedge\sup\L_U}}
\sum_{t_1\leq s_1\in\L_H^{\leq s}}\mu_H(t_1,s_1)\mu(U,H)\\
\vspace{0.1cm}
\mbox{}&\displaystyle
=\left\{
\begin{array}{cl}
\vspace{0.1cm}\displaystyle
\sum_{U\leq H\leq K}\mu(U,H)\quad&\mbox{if }t=s\wedge\sup\L_U,\\
\vspace{0.1cm}
0\quad&\mbox{otherwise}
\end{array}\right. \\ 
\mbox{}&\displaystyle
=\left\{
\begin{array}{cl}
\vspace{0.1cm}\displaystyle
1\quad&\mbox{if }(U,t)=(K,s),\\
0\quad&\mbox{otherwise}
\end{array}\right.
\end{array}
\]
for all $(U,t)\in\XX$. If $M_\L=M_\SSS$ (see Example \ref{epl:slice}), then 
this equality is stated in \cite[Proposition 5.3]{Bou12}, and the primitive 
idempotent $e_{(K,s)}^{(G)}$ is due to \cite[Corollary 5.5]{Bou12}. 
\end{remark}
\par
Assume that $(K\cap U,s\wedge t)\in\XX$ for all $(K,s),\,(U,t)\in\XX$ and 
$(G,\sup\L_G)\in\XX$. By Proposition \ref{pro:monoid} and 
Theorem \ref{thm:partial p}, $\mathrm{\Omega}(G,\XX)_{(p)}$ is a subring of 
$\mathrm{\Omega}(G,M_\L)_{(p)}$, which is called the partial lattice Burnside 
ring relative to $\XX$. 
\par 
We identify $\mathrm{\Omega}(G,\XX)$ with $\bz\boldGamma_\XX$, and turn to 
a criteria for the units of $\mathrm{\Omega}(G,\XX)$. 
\begin{proposition}
\label{pro:Yoshida criteria partial}
Assume that $(K\cap U,s\wedge t)\in\XX$ for all $(K,s),\,(U,t)\in\XX$ and 
$(G,\sup\L_G)\in\XX$. Let $\widetilde{x}=(x_{(U,t)})_{(U,t)\in\ov{\XX}}\in
\mathrm{\Omega}(G,\XX)^\times$. Then $\widetilde{x}$ is contained in the image 
of the ring homomorphism $\varphi_\XX:\mathrm{\Omega}(G,\XX)\to
\widetilde{\mathrm{\Omega}}(G,\XX)$, where 
$\varphi_\XX=\varphi_\XX^{(\infty)}$, if and only if, for each 
$(U,t)\in\ov{\XX}$, the map 
$\gamma_{(U,t)}^{\widetilde{x}}:\aut((G/U)_t)\to\lr{-1}$ given by 
\[
\sigma_g\mapsto x_{(U,t)}x_{(\ov{\lr{g}U},\ov{s_{(g,t)}})} 
\]
for all $gU\in W_G(U,t)$ is a group homomorphism. 
\end{proposition}
{\itshape Proof.}
Since $\mathrm{\Omega}(G,\XX)$ is a subring of $\mathrm{\Omega}(G,M_\L)$, 
the proposition is a consequence of Proposition \ref{pro:Yoshida}, 
Corollary \ref{cor:ring homomorphism}, and Lemma \ref{lem:partial}. 
$\Box$
\begin{example}
\label{epl:section unit}
We set $\XX=\{(K,E)\st K\unlhd E\leq G\}$, where $K\unlhd E$ denotes that $K$ 
is a normal subgroup of $E$. Then $\XX$ is a subset of $\SS(G,M_\SSS)$ closed 
under the action of $G$, where $M_\SSS$ is given in Example \ref{epl:slice}. 
Let $(U,F)\in\XX$, and let $gU\in W_G(U,F)$. We denote by 
$(\lr{g}U)^{\unlhd\lr{g}F}$ the normal closure of $\lr{g}U$ in $\lr{g}F$. If 
$(K,E)\geq(\lr{g}U,s_{(g,F)})$ with $(K,E)\in\XX$, then 
$E\geq s_{(g,F)}=\lr{g}F$ and $K\geq(\lr{g}U)^{\unlhd\lr{g}F}\geq\lr{g}U$, 
which implies that 
$(\ov{\lr{g}U},\ov{s_{(g,F)}})=((\lr{g}U)^{\unlhd\lr{g}F},\lr{g}F)$. Hence 
the condition (A) holds. The partial lattice Burnside ring 
$\mathrm{\Omega}(G,\XX)$ is isomorphic to the section Burnside ring 
$\mathrm{\Gamma}(G)$ introduced by S. Bouc \cite{Bou12}. The primitive 
idempotents of $\bq\otimes_\bz\mathrm{\Gamma}(G)$, together with 
the description of $\mu_{\XX}$ in terms of the M\"{o}bius function of 
$\SSS(G)$, are given in \cite[Corollary 12.5]{Bou12}. Obviously, 
$(K\cap U,E\cap F)\in\XX$ for all $(K,E),\,(U,F)\in\XX$ and 
$(G,G)\in\XX$, so that Proposition \ref{pro:Yoshida criteria partial} in this 
case is equivalent to \cite[Theorem 15.3]{Bou12}. 
\end{example}
%----------------------------------------------------------------
%
%    References 
%
%------------------------------------------------------- 
 
%------------------------------------------------------- 
%    END OF REFERENCES 
%------------------------------------------------------- 
%\input{latticex}
\end{document}